\documentclass[12pt, twoside, reqno]{amsart}
\usepackage{amssymb}
\usepackage{amsmath}
\usepackage{amsthm}
\usepackage{graphicx, subcaption}
\usepackage{float, cases}
\usepackage[active]{srcltx}
\usepackage{forest}
\usepackage[hidelinks]{hyperref}
\usepackage{amsfonts}
\usepackage{amssymb}
\usepackage{amsmath,bbm,amssymb,amsxtra}
\usepackage{mathrsfs}
\usepackage{enumerate}
\usepackage{caption}
\allowdisplaybreaks
\usepackage{epsfig}
\usepackage{ae, mdframed}

\newtheorem{theorem}{Theorem}[section]

\newtheorem{lemma}[theorem]{Lemma}

\newtheorem{definition}[theorem]{Definition}
\newtheorem{example}[theorem]{Example}

\newtheorem{theorema}{Theorem}

\newtheorem{problema}{Problem}

\newtheorem{lemmaa}{Lemma}

\newcommand{\IM}{\text{Im}}

\newcommand{\diam}{\text{diam}}

\newcommand{\ra}{\rightarrow}
\newcommand{\mbs}{\mathbb{S}}
\newcommand{\mbd}{\mathbb D}

\newcommand{\mbn}{\mathbb N}
\newcommand{\mcc}{\mathcal C}
\newcommand{\mbh}{\mathbb H}
\newcommand{\mbr}{\mathbb{R}}
\newcommand{\pao}{\partial \Omega}

\numberwithin{equation}{section} 
\makeatletter
\@addtoreset{equation}{section}
\makeatother

\makeatletter
\renewenvironment{proof}[1][\proofname] {\par\pushQED{\qed}\normalfont\topsep6\p@\@plus6\p@\relax\trivlist\item[\hskip\labelsep\bfseries#1\@addpunct{.}]\ignorespaces}{\popQED\endtrivlist\@endpefalse}
\makeatother

\makeatletter
\@namedef{subjclassname@2020}{%
  \textup{2020} Mathematics Subject Classification}
\makeatother

\makeatletter
\newcommand*{\mint}[1]{%
  \mint@l{#1}{}%
}
\newcommand*{\mint@l}[2]{%
  \@ifnextchar\limits{%
    \mint@l{#1}%
  }{%
    \@ifnextchar\nolimits{%
      \mint@l{#1}%
    }{%
      \@ifnextchar\displaylimits{%
        \mint@l{#1}%
      }{%
        \mint@s{#2}{#1}%
      }%
    }%
  }%
}
\newcommand*{\mint@s}[2]{%
  \@ifnextchar_{%
    \mint@sub{#1}{#2}%
  }{%
    \@ifnextchar^{%
      \mint@sup{#1}{#2}%
    }{%
      \mint@{#1}{#2}{}{}%
    }%
  }%
}
\def\mint@sub#1#2_#3{%
  \@ifnextchar^{%
    \mint@sub@sup{#1}{#2}{#3}%
  }{%
    \mint@{#1}{#2}{#3}{}%
  }%
}
\def\mint@sup#1#2^#3{%
  \@ifnextchar_{%
    \mint@sup@sub{#1}{#2}{#3}%
  }{%
    \mint@{#1}{#2}{}{#3}%
  }%
}
\def\mint@sub@sup#1#2#3^#4{%
  \mint@{#1}{#2}{#3}{#4}%
}
\def\mint@sup@sub#1#2#3_#4{%
  \mint@{#1}{#2}{#4}{#3}%
}
\newcommand*{\mint@}[4]{%
  \mathop{}%
  \mkern-\thinmuskip
  \mathchoice{%
    \mint@@{#1}{#2}{#3}{#4}%
        \displaystyle\textstyle\scriptstyle
  }{%
    \mint@@{#1}{#2}{#3}{#4}%
        \textstyle\scriptstyle\scriptstyle
  }{%
    \mint@@{#1}{#2}{#3}{#4}%
        \scriptstyle\scriptscriptstyle\scriptscriptstyle
  }{%
    \mint@@{#1}{#2}{#3}{#4}%
        \scriptscriptstyle\scriptscriptstyle\scriptscriptstyle
  }%
  \mkern-\thinmuskip
  \int#1%
  \ifx\\#3\\\else_{#3}\fi
  \ifx\\#4\\\else^{#4}\fi
}
\newcommand*{\mint@@}[7]{%
  \begingroup
    \sbox0{$#5\int\m@th$}%
    \sbox2{$#5\int_{}\m@th$}%
    \dimen2=\wd0 %
    \let\mint@limits=#1\relax
    \ifx\mint@limits\relax
      \sbox4{$#5\int_{\kern1sp}^{\kern1sp}\m@th$}%
      \ifdim\wd4>\wd2 %
        \let\mint@limits=\nolimits
      \else
        \let\mint@limits=\limits
      \fi
    \fi
    \ifx\mint@limits\displaylimits
      \ifx#5\displaystyle
        \let\mint@limits=\limits
      \fi
    \fi
    \ifx\mint@limits\limits
      \sbox0{$#7#3\m@th$}%
      \sbox2{$#7#4\m@th$}%
      \ifdim\wd0>\dimen2 %
        \dimen2=\wd0 %
      \fi
      \ifdim\wd2>\dimen2 %
        \dimen2=\wd2 %
      \fi
    \fi
    \rlap{%
      $#5%
        \vcenter{%
          \hbox to\dimen2{%
            \hss
            $#6{#2}\m@th$%
            \hss
          }%
        }%
      $%
    }%
  \endgroup
}

\begin{document}

\title{Sobolev homeomorphisms and Jordan curves with positive area}

\author[P. Koskela]{Pekka Koskela}
\address{P. Koskela, Department of Mathematics and Statistics, University of Jyväskylä, P.O. Box 35 (MaD), FI-40014, Finland}
\email{pekka.j.koskela@jyu.fi}

\author[H. Xu]{Haiqing Xu}
\address{H. Xu, Frontiers Science Center for Nonlinear Expectations (Ministry of Education of China), Research Center for Mathematics and Interdisciplinary Sciences, Shandong University, Qingdao 266237, China}
\email{haiqing.xu@email.sdu.edu.cn}

\thanks{P.  K.  was partially supported by the Finnish centre of excellence in Randomness and Structures (FiRST), Research Council of Finland project number 364210.
H. X. was supported by the Natural Science
Foundation of Shandong Province (No.~ZR2025MS96),  and by
the Fundamental Research Funds for the Central Universities.
}

\begin{abstract}
Each homeomorphic parametrization of a Jordan curve via the
unit circle $\mbs\subset \mbs^2$ extends to a homeomorphism of the entire $\mbs^2.$
We show that there exists a  Jordan curve of positive area whose every parametrization via $\mbs$ admits a homeomorphic extension that has $p$-Sobolev regularity on $\mbs^2$ for all exponents $1\le p<2 .$
This is obtained via a new sharp extension result for the case of $p<2.$
\end{abstract}
\maketitle

\section{Introduction}

In this paper, we consider Jordan curves in the plane $\mbr^2 $ and on the (Riemann) sphere $\mbs^2.$ They are homeomorphically parametrized by the unit circle $\mbs.$
By the Jordan-Sch\"onflies theorem each such a parametrization extends to a self-homeomorphism of $\mbr^2 ,$ or of   $\mbs^2$  when we consider  $\mbs$  embedded in $\mbs^2.$ 
The extensions given by this theorem do not need to have any regularity even when both the parametrization and the curve are regular.
In this paper, we focus on a Sobolev variant of the Jordan–Sch\"onflies theorem.

\begin{problema}\label{25.3.19_11}
For which Jordan curves in  $\mbs^2$  and $1\le p\le \infty,$ does every parametrization admit a homeomorphic extension  
$h:   \mbs^2 \xrightarrow{\mathrm{onto}} \mbs^2$ in the class $W^{1,p} (\mbs^2 ,\mbs^2)?$
\end{problema}

The answer necessarily depends on the Jordan curve in question even when $p=1$ since there exist examples that do not admit any such Sobolev extensions \cite{Yi}. 
Also, each homeomorphism with $W^{1,2}-$regularity satisfies Luzin's condition (N): sets of area zero get mapped to sets of area zero \cite{Hencl 2014, Resetnjak 1966}. Hence, no parametrization of a Jordan curve of positive area can have a homeomorphic $W^{1,2}$-extension. In fact, the above question is relevant (even for parametrizations of $\mbs$) only when $1\le p<2,$ since each homeomorphism with $W^{1,2}-$regularity has a logarithmic modulus of continuity
\cite{Hencl 2014, Reshetnyak 1989}.

Our first result shows that, surprisingly, there exist Jordan curves of positive area that are as good for this problem as the unit circle $\mbs$ is.

\begin{theorem}\label{posa}
There exists a Jordan curve of positive area in $\mbs^2$ whose every parametrization via $\mbs \subset \mbs^2$ has a homeomorphic extension that belongs to $W^{1,p} (\mbs^2 ,\mbs^2)$ for each $1\le p<2.$
\end{theorem}

One should compare this existence result to an old construction of a Sobolev homeomorphism with $W^{1,p}$-regularity, $1\le p<2,$ that fails to satisfy Luzin's condition (N) \cite{ponomarev}. In that construction, a product Cantor set gets mapped onto a Cantor set of positive area.  In our case, the Jordan curve in question contains a product-type Cantor set of positive area and its preimage is then contained in 
$\mbs.$

The proof of Theorem \ref{posa} is based on our sharp extendability results in terms of $p$ and the growth of the (quasi)hyperbolic distance. Given a planar Jordan curve, we call the bounded complementary component a Jordan domain. Then its boundary is the given Jordan curve.
In what follows, $\mbd$ is the open unit disk and $\rho_{\Omega}(z,z_0)$ is the hyperbolic distance between points $z,z_0$ in a Jordan domain $\Omega.$

\begin{theorem}\label{perus2}
 Let $\Omega$ be a Jordan domain. Suppose that for every $0<\alpha<1$  there is a constant $C_{\alpha}$ so that the hyperbolic distance $\rho_{\Omega}$  satisfies the growth condition
$$\rho_{\Omega}(z,z_0)\le C_{\alpha} d(z,\partial  \Omega)^{-\alpha}$$ 
for some fixed $z_0\in \Omega$ and all $z\in \Omega.$ Given a homeomorphism $\varphi:\mbs\to \partial \Omega,$ there is a homeomorphic extension $f:\overline{\mbd}\to \overline{\Omega}$ that belongs to
$W^{1,p}(\mbd,\Omega)$ for all $1\le p<2.$ 
\end{theorem}

Our approach also allows us to give an optimal sufficient conditions for homeomorphic $W^{1,p}-$extendability in terms of the growth of the hyperbolic metric.

\begin{theorem}\label{perus}
 Let $\Omega$ be a Jordan domain and $0<\alpha<1.$ Suppose that there is a constant $C$ so that the hyperbolic distance $\rho_{\Omega}$  satisfies the growth condition
$$\rho_{\Omega}(z,z_0)\le C d(z,\partial  \Omega)^{-\alpha}$$ 
for some fixed $z_0\in \Omega$ and all $z\in \Omega.$ Given a homeomorphism $\varphi:\mbs\to \partial \Omega,$ there is a homeomorphic extension $f:\overline{\mbd}\to \overline{\Omega}$ that belongs to
$W^{1,p}(\mbd,\Omega)$ for all $1\le p<2-\alpha.$ Moreover, this range of $p$ is sharp.
\end{theorem}

To show the sharpness,  for each $\alpha \in (0,1)$ we construct a Jordan domain $\Omega_{\alpha}$ with the above geometric assumption and one parametrization of $\partial \Omega_{\alpha}$ via $\mbs$ such that the Sobolev regularity of its homeomorphic extensions fails to reach the threshold $2-\alpha .$

In analysis, self-improving properties are typical: if one has uniform estimates for a function class, then something slightly better also holds. We modify our construction so as to show that this is not the case for 
homeomorphic extensions.

\begin{theorem}\label{26.7.20}
Given $1<p<2,$ there is a Jordan curve $\gamma_p$ in $\mbs^2$ and a constant $M$ so that every parametrization of $\gamma_p$ via $\mbs \subset \mbs^2$ has an extension to a self-homeomorphism of  $\mbs^2$ with $W^{1,p}-$norm bounded by $M, $ but there is a parametrization that does not have a $W^{1,q}-$extension for any $q>p.$
\end{theorem}

The study of homeomorphic Sobolev extensions has been active in recent years \cite{Koskela 2020, Koskela 2018,Koski 2021, Koski 2023, Koski 2025, Wang 2023, Xu 2020}.  The result closest to Theorem \ref{perus} is that the $q$-integrability, $q>1,$ of the hyperbolic distance to a base point $z_0$ guarantees the existence of a homeomorphic extension to $\mbd$ with $W^{1,p}$-regularity for all $1\le p<2,$ see \cite{Bouchala 2025}. We cannot apply this result since the optimal pointwise growth condition guaranteeing this integrability condition is $h_{\Omega}(z,z_0)\le C \left( \log (d(z_0 , \pao) / d(z ,\pao) ) \right)^2 ,$ see  \cite{Koskela 2005}. This latter growth condition actually implies that the area of $\partial \Omega$ must be zero, and so does the above $q$-integrability condition on $h_{\Omega}.$  This is not necessarily the case for a domain $\Omega$ in Theorem \ref{perus} that gives the first known examples where the extendability property only holds for exponents $p$ close to 1. 
 
 The positive direction in Theorem \ref{perus} is obtained by showing that our assumptions allow us to employ the general extension result from \cite{Koski 2023}.
 The examples towards sharpness are based on modifications to the classical construction of a Jordan curve of positive area from \cite{Osgood 1903}. Theorem \ref{posa} is then a consequence of Theorem \ref{perus} applied to an even more novel construction of a Jordan curve of positive area.
 
 Partial motivation for this work comes from the variational approach to geometric function theory \cite{Astala 2009, Iwaniec 2001} and
nonlinear elasticity \cite{Ball 1976}. One aims to find an energy minimizer for a boundary value problem with homeomorphic boundary values. The competing class consists of homeomorphic Sobolev maps 
of finite energy and with the given homeomorphic boundary value. In order for this to make sense, one needs to know that the competing class is non-empty. Theorem \ref{posa} shows that the minimization problem makes
sense even for arbitrary parametrizations of certain Jordan curves of positive area.

This article is structured as follows.
In Section \ref{25.3.19_1}, we collect the preliminaries required to prove our theorems. They include the basic properties of the hyperbolic metric and the quasihyperbolic metric.
In Section \ref{25.4.1_1}, we prove Theorem \ref{perus} and Theorem \ref{perus2}. The proof of Theorem \ref{26.7.20} is given in Section \ref{26.8.11_1}. Finally, Section \ref{26.4.10_1} is dedicated to the proof of Theorem \ref{posa}.

\section{Preliminaries}\label{25.3.19_1}


This section contains the necessary definitions and some basic results.

\begin{definition}
Assume that $\Omega$ is an open subset of $\mathbb{R}^2 .$
The Sobolev space $W^{1,p}(\Omega)$ for $1\le p <\infty$ consists of all $p$-integrable functions $u$  that have $p$-integrable first order distributional derivatives on $\Omega.$ In symbols,
$$W^{1,p}(\Omega):=\left\{ u \in L^p(\Omega): |D u |\in L^p(\Omega) \ \right\}.$$	
\end{definition}

Continuous functions in $W^{1,p} (\Omega)$ are conveniently characterized in terms of absolute continuity on lines.
For $\nu=1$ and $2 ,$ let $I_\nu=[a_\nu,b_\nu]$ be closed intervals such that the associated rectangle $Q=I_1\times I_2$ lies in $\Omega .$ Consider the faces
$$Q_1=\{a_1 \}\times I_{2} \qquad \mbox{ and } \qquad  Q_2= I_{1} \times \{a_2 \} .$$
We say that a continuous function $f:\Omega\to\mathbb{R}^m$ has the  $ACL-$property  {\it  or is absolutely continuous on lines}, if for each such $Q,$ both $\nu = 1$ and $2$ and almost every $a\in Q_\nu$, with respect to $\mathcal{L}^1$-measure, the function
$$t\mapsto f(a+te_\nu),\quad0\leqslant t\leqslant b_\nu-a_\nu $$
is absolutely continuous. Here, of course, $e_\nu$ denotes the $\nu^\mathrm{th}$ unit basis vector.
Clearly, such functions have classical first-order partial derivatives almost everywhere in $\Omega$ and these derivatives are Borel functions.  The following result can be found e.g., in  \cite{Rickman 1993}.

\begin{theorema}\label{1-thm 4.1.10}
Let $u \in L^p(\Omega;\mathbb R^2)$ be continuous, where  $1\le p <\infty .$ Then $u$ lies in $W^{1,p}(\Omega;\mathbb R^2)$ if and only if $u$ has the ACL-property and the (classical) partial derivatives belong to $L^p(\Omega).$ 
\end{theorema}


The natural variant of the above also holds for the case $\Omega \subset \mbs^2$ as seen by reducing it to the planar case via bi-Lipschitz coordinate patches.

 We will also need the following fact. If $u\in W^{1,p}(\Omega;\mathbb R^2)$ is continuous, then $u$ restricted to $S^1(x,r)$ is absolutely continuous for each $x\in \Omega$ and almost every $0<r<d(x,\partial \Omega)$
 
Next, we recall the definitions of the hyperbolic and quasihyperbolic metrics and a couple of basic properties of conformal mappings.

\begin{definition}
The hyperbolic metric $\rho_{\mbd}$ on the unit disk $ \mbd $ is defined by 
\begin{equation*}\label{0112_1}
\rho_{\mbd} (z_1 ,z_2) = \inf_{\gamma} \int_{\gamma} \frac{2 \, |dz|}{1-|z|^2} ,
\end{equation*}
where the infimum is taken over all rectifiable arcs $\gamma$ joining $z_1$ to $z_2$ in $\mbd .$
\end{definition}
The infimum is attained by the {\it hyperbolic geodesic} from $z_1$ to $z_2 ,$ that is the image of the non-Euclidean segment from $z_1$ to $z_2 .$
There is a close relationship between hyperbolic metrics and conformal mappings.
\begin{lemmaa}{{\rm (\cite[Theorem 2.10.8]{Astala 2009})}}\label{23.12.7_3}
Suppose that $f$ is conformal in the unit disk $\mbd$ and $z ,w \in \mbd .$ Then 
\begin{equation*}
\exp (-3 \rho_{\mbd} (z,w)) \le \frac{|f'(z)|}{|f'(w)|} \le \exp (3 \rho_{\mbd} (z,w)).
\end{equation*}
\end{lemmaa}
The above allows one to consider the derivative of a conformal mapping to be almost constant on hyperbolic disks.

Let $\Omega$ be a simply connected domain. By the Riemann mapping theorem, there is a conformal mapping $f$ from $\mbd$ onto $\Omega .$ The hyperbolic metric $\rho_{\Omega}$ on $\Omega$ is defined as 
\begin{equation*}\label{0112_3}
\rho_{\Omega} (z_1 ,z_2) := \rho_{\mbd} (f^{-1} (z_1), f^{-1} (z_2)).
\end{equation*} 
Then $\rho_{\Omega}$ is well-defined because $\rho_{\mbd}$ is invariant under  M\"obius transformations.
It also follows that the hyperbolic metric is conformally invariant.
Usually, one does not have an explicit formula for the above $f ,$ and thus neither for $\rho_{\Omega}.$
Next, we introduce a useful substitute for $\rho_{\Omega} .$ 

\begin{definition}
The quasihyperbolic metric $h_{\Omega}$ on a domain $\Omega \subset \mbr^2$ is defined as 
\begin{equation*}\label{0113_5}
h_{\Omega} (z_1 ,z_2) = \inf_{\gamma} \int_{\gamma} \frac{1}{d(z, \pao)} \, |dz|,
\end{equation*}
where the infimum is taken over all rectifiable arcs $\gamma$ joining $z_1$ to $z_2$ in $\Omega ,$ and $d(z, \pao)$ is the Euclidean distance between $z$ and $\pao .$
\end{definition}
 It follows from \cite[Lemma 2.1]{Gehring 1976} that 
\begin{equation}\label{25.3.11_1}
h_{\Omega} (z_1 ,z_2) \ge \log \left( \max\left\{ \frac{d(z_1 ,\pao)}{d(z_2 ,\pao)},\ \frac{d(z_2 ,\pao)}{d(z_1 ,\pao)} \right\} \right)
\end{equation}
for all $z_1 ,\ z_2 \in \Omega .$

We introduce the following form of Koebe's $1/4$-theorem. 
\begin{theorema}{{\rm (\cite[Theorem 2.10.6]{Astala 2009})}}\label{23.12.7_8}
Suppose that $f$ is conformal in a domain $\Omega$ with $f(\Omega) = \Omega' \subset \mbr^2 .$ Let $z \in \Omega .$ Then 
\begin{equation*}
\frac{1}{4} |f' (z)| d (z ,\pao) \le d(f(z), \pao') \le |f' (z)| d (z ,\pao).
\end{equation*}
\end{theorema}
By this theorem, we see that $h_{\Omega} (z_1 ,z_2) \approx \rho_{\Omega} (z_1 ,z_2) $ for all $z_1 ,\ z_2 \in \Omega .$

Let us recall the definition of $N-$functions.
A real-valued function $\Phi$ defined on $[0,\infty)$ is called an {\it $N-$function}, if it admits the representation 
\begin{equation*}
\Phi (t) = \int^t _0 \varphi (\tau) d \tau,
\end{equation*}
where the function $\varphi (\tau)$ is right-continuous for all $\tau \ge 0 ,$  positive for $\tau >0,$ non-decreasing and satisfies the conditions $\varphi (0) =0 $ and $ \lim_{\tau \ra \infty} \varphi (\tau) =\infty .$
Given such a function $\varphi ,$
we define the function $\psi (\xi)$ for $\xi \ge 0$ by the equality
\begin{equation}\label{26.7.20_2}
\psi (\xi) = \sup_{\varphi (\tau) \le \xi} \tau .
\end{equation}
It is easily seen that the function $\psi$ also satisfies the above properties of $\varphi .$ 
Denote 
$$\Psi (s) = \int^s _{0} \psi (\xi) d \xi  .$$ 
Then $\Phi$ and $\Psi$ are called {\it mutually complementary $N-$functions}.
Examples of such complementary pairs are $\Phi(t) = t^{p} /p$ and $\Psi (s) = s^{p'} /p' ,$ where $p \in (1 ,\infty)$ and $p'$ is the conjugate index of $ p.$

Let $\Phi$ be an $N-$function. The {\it Orlicz class} $L^{\dagger}_{\Phi}$ is the set of those real-valued measurable functions, defined on $\mbr^n ,$ for which $\int_{\mbr^n} \Phi (|f(x)|) dx <\infty .$ The $L^{\dagger}_{\Phi}$  is not necessarily a vector space, because it might fail to be closed under the scalar multiplication.  The {\it Orlicz space} $L_{\Phi}$  is the linear hull of the Orlicz class  $L^{\dagger}_{\Phi}, $ that is, the smallest linear space containing $L^{\dagger}_{\Phi} .$
The reader may verify that 
\begin{equation*}
||f||_{\Phi} = \inf \left\{k >0 : \int_{\mbr^n} \Phi \left(\frac{|f(x)| }{k}\right) dx \le 1 \right\}
\end{equation*}
is a norm on $L_{\Phi} .$ 
It is called the {\it Luxemburg norm}. The infimum is attained, and $ \int_{\mbr^n} \Phi (|f (x)| /||f||_{\Phi}) = 1 .$
For this theory, see e.g.  \cite[Chapter II]{Krasnoselskii 1961}.

Let $\Phi$ and $\Psi$ be mutually complementary $N-$functions.
It follows from \eqref{26.7.20_2} that 
\begin{equation*}\label{26.7.21_1}
\psi (\varphi (\tau)) \ge \tau \qquad \mbox{ and } \qquad\varphi(\psi(\xi)) \ge \xi 
\end{equation*}
for any $\tau \ge 0$ and any $\xi \ge 0 .$
If additionally the function $\varphi$ is continuous and strictly increasing, then the function $\psi$ is the ordinary inverse function of $\varphi .$
Consider the Cartesian coordinate system with the horizontal coordinate $t$ and the vertical coordinate $s .$
Then the value $ts$ represents the area of a rectangle with side lengths $t$ and $s ,$ and $\Phi (t)$ is the area of lower part of the graph $\xi = \varphi (\tau) $ from $\tau =0$ to $\tau =t .$ Analogously, $\Psi (s)$ is the area of the left part of the graph $\tau = \psi(\xi)$ from $\xi= 0$ to $\xi =s .$
It is geometrically clear that 
\begin{equation*}
ts \le \Phi (t) + \Psi (s) 
\end{equation*}
for all $t \ge 0$ and $s \ge 0.$ This is called  {\it Young's inequality}.
Apply Young's inequality to $ f(x) / ||f||_{\Phi} $ and $g(x) /  ||g||_{\Psi} $ and then integrate over $\mbr^n .$
We obtain
\begin{equation}\label{26.7.20_3}
\left | \int_{\mbr^n} f(x) g(x) dx\right | \le 2 ||f||_{\Phi} ||g||_{\Psi} .
\end{equation}
This is called a {\it generalized H\"older's inequality}.
One can find this inequality, e.g., in  \cite[Chapter 8]{Adams 2003}.

In this paper, we need a discrete variant of the generalized H\"older's inequality.
Let $a= \{a_i\}^{\infty} _{i=1}$ and $b=\{b_i\}^{\infty} _{i=1}$ be two sequences of numbers.  Decompose $\mbr^n$ into countable number of cubes $\{Q_i\}^{\infty} _{i=1}$ of unit size. 
Apply the inequality \eqref{26.7.20_3} to $f(x) = \sum^{\infty} _{i=1} a_i \chi_{Q_i} (x)$ and $g(x) = \sum^{\infty} _{i=1} b_i \chi_{Q_i} (x) .$
It follows that 
\begin{equation}\label{26.7.21_4}
\left|\sum^{\infty}_{i=1}  a_i b_i\right| \le 2 ||a||_{\Phi} ||b||_{\Psi},
\end{equation}
where
$||a||_{\Phi} =  \inf \{k >0 : \sum^{\infty}_{i=1}\Phi (|a_i| /k ) \le 1 \}$
and 
$||b||_{\Psi} =  \inf \{k >0 : \sum^{\infty}_{i=1}\Psi ( |b_i| /k ) \le 1 \} .$

\section{Proofs of  Theorem \ref{perus} and Theorem \ref{perus2}}\label{25.4.1_1}

We split the proof of Theorem \ref{perus} into two parts.
The proof of the positive part relies on a sufficient condition for homeomorphic Sobolev extensions from \cite{Koski 2023}. We reformulate this result in a form that is convenient for us.
In what follows, $\ell(\gamma)$ refers to the length of a curve $\gamma.$

\begin{theorema}{{\rm (\cite[Theorem 3.1]{Koski 2023})}}\label{23.12.10_1}
Suppose that $\Omega$ is a Jordan domain and $\varphi : \partial \mbd \ra \partial \Omega$ is a boundary homeomorphism. 
Suppose that there is a dyadic family $I = \{I_{n,j} \},$ $n=1,2,\dots,$ $j=1,\dots, 2^n,$ of closed
arcs on $\partial \mbd$ and associated crosscuts $\Gamma_{n,j}$ in $\Omega$ connecting the two endpoints of the boundary arc $\varphi(I_{n,j} ) \subset \partial \Omega$ so that the crosscuts are pairwise disjoint apart from their endpoints at the boundary.  Suppose additionally that $\lim_{n\to \infty}\sup_j \ell(\Gamma_{n,j})=0.$ Then $\varphi$ admits a homeomorphic extension from $\overline {\mbd}$ to $\overline \Omega$ that belongs to the class $W^{1,p}(\mbd,\Omega)$ for all those $1<p<\infty$ for which
$$
\sum^{\infty}_{n=n_0} 2^{(p-2)n} \sum^{2^n}_{j=1} \ell(\Gamma_{n,j})^p <\infty
$$
holds for some $n_0\in N.$
Moreover, the Sobolev norm is bounded by the double sum.
\end{theorema}

\begin{proof}[Proof of  the positive part of Theorem \ref{perus}]
Fix a dyadic family $I$ of closed arcs on $\partial \mbd$. In order to apply Theorem \ref{23.12.10_1}, we construct
crosscuts $\{\Gamma_{n,j}\}$ such that $\lim_{n\to \infty}\sup_j \ell(\Gamma_{n,j})=0$ and 
$\sum^{\infty}_{n=n_0} 2^{(p-2)n} \sum^{2^n}_{j=1} \ell(\Gamma_{n,j})^p <\infty$ for all $p<2-\alpha .$
Let $\Gamma_{n,j}$ be the hyperbolic geodesic connecting the endpoints of $\varphi (I_{n,j}) .$
It remains  to estimate the lengths  $\ell(\Gamma_{n,j}) . $

Fix a conformal mapping $f :\mbd \ra \Omega$ with $f(0) =z_0 .$ 
Denote by  $\xi_1 ,\ \xi_2 \in \mbs$ the preimages of the endpoints of $\varphi (I_{n,j})$ under $f .$
We may assume that $\xi_1 = e^{i \pi /4}$ and $\xi_2 = e^{-i \pi /4} .$
Let $\mbh$ be the upper-half plane.
Define a conformal mapping $T :\mbd \ra \mbh$ as $T(z) = a (1-z)/(1+z) $ with $a = (1+\xi_1) /(1- \xi_1) .$
Denote 
\begin{equation*}
g = f \circ T^{-1} .
\end{equation*}
Then $g : \mbh \ra \Omega$ is a conformal homeomorphism.
Since $T(\xi_1) = 1=-T(\xi_2),$ the mapping $T$ maps the hyperbolic geodesic joining $\xi_1$ to $\xi_2$ onto the half circle $\mbs \cap \mbh .$
Let 
$$\gamma_m = \{e^{i \theta} : \theta \in (\theta_{m+1} , \theta_m ]\} $$
with $\theta_m = \pi 2^{-m}$ for $m \ge 1  .$
Then $\cup^{\infty} _{m=1} \gamma_m$ is half of the geodesic connecting $T(\xi_1) $ and $T(\xi_2) .$
Thus $\cup^{\infty} _{m=1} g(\gamma_m) $ is a part of $\Gamma_{n,j} .$
We next estimate the length $\ell(g(\gamma_m))$ and then $\ell(\cup^{\infty} _{m=1} g(\gamma_m) ) .$ Analogously, there will be a length estimate for $\Gamma_{n,j} \setminus  \cup^{\infty} _{m=1} g(\gamma_m) .$ The combination of these will give an estimate for $\ell(\Gamma_{n,j}). $

Let  $r_m = 1-2^{-m-1} $ and 
\begin{equation*}
A_m = \{r e^{i \theta} : r \in [r_m ,1] \mbox{ and } \theta \in (\theta_{m+1} , \theta_m ]\} 
\end{equation*}
for $m \ge 1 .$
Set $z_m = e^{i \theta_m} .$ 
For $z \in A_m ,$ we denote by $L_{z_m ,z}$ the line segment connecting $z_m$ and $z .$
Then  the quasihyperbolic distance $h_{\mbh}$ in $\mbh$ satisfies
\begin{equation}\label{23.12.7_4}
h_{\mbh} (z_m ,z) \le  \int_{L_{z_m ,z}} \frac{1}{\IM \, z} |dz| \le   \frac{\ell(L_{z_m ,z})}{r_m \sin \theta_{m+1}} \le C ,
\end{equation}
where the constant $C$ is independent of $m .$
Recall that $h_{\mbh}$ is comparable to the hyperbolic metric.
It follows from Lemma \ref{23.12.7_3} that $|f' (T^{-1} (z_m))| \approx |f' (T^{-1} (z))|$ for all $z \in A_m .$
Notice that $|(T^{-1})' (z)| \approx |(T^{-1})' (z_m)| \approx 1/|a|$ for all $z \in A_m .$ It follows that $|g'(z_m)| \approx  |g'(z)|$
for all $z \in A_m .$
Hence, 
\begin{equation}\label{23.12.7_5}
\ell(g(\gamma_m)) =  \int_{\gamma_m} |g'(z)| \, |dz| \approx  |g'(z_m)| 2^{-m} 
\end{equation}
and 
\begin{equation}\label{23.12.7_6}
|g(A_m)| = \int_{A_m} |g'(z)|^2 dz \approx |g'(z_m)|^2 2^{-2m} ,
\end{equation}
where $|A_m|$ is the area of $A_m.$
Since $\IM\, z \approx 2^{-m} $ for $z \in A_m ,$ by \eqref{25.3.11_1} we obtain that $h_{\mbh} (T(0),z) \gtrsim m .$
Notice that the quasihyperbolic metric is preserved under conformal mappings modulo multiplicative constants. It follows that 
$h_{\Omega} (z_0, g(z))\gtrsim m $
for all $z \in A_m .$
Fix $q \in (2/(2-\alpha) , 2/\alpha).$
We obtain from \eqref{23.12.7_5} and \eqref{23.12.7_6} that 
\begin{equation}\label{25.3.11_2}
\ell(g(\gamma_m))^2 m^q  \approx  |g(A_m)|m^q \lesssim \int_{g(A_m)} h^q _{\Omega} (z_0 ,\xi) d \xi .
\end{equation}
Moreover, by our assumption on the growth of the (quasi)hyperbolic metric and Theorem \ref{23.12.7_8}, we have
\begin{align}
\int_{g(A_m)} h^q _{\Omega} (z_0 ,\xi) d \xi
\lesssim & \int_{g(A_m)} \frac{1}{d^{\alpha q} (\xi ,\pao)} d \xi
\approx \int_{g(A_m)} \frac{|(g^{-1})' (\xi)|^{\alpha q }}{d^{\alpha q }(g^{-1} (\xi) ,\partial \mbh) } d \xi \notag  \\
= & \int_{A_m} \frac{|g' (z)|^{2-\alpha q }}{d^{\alpha q }(z ,\partial \mbh) } d z \notag \\
\approx & \frac{|g' (z_m)|^{2-\alpha q }}{2^{-m(\alpha q -2)}} \approx |g(A_m)| ^{1-\frac{\alpha q}{2}} ,\label{25.3.11_3}
\end{align}
where the equality comes from a change of variable.
Therefore, combining \eqref{25.3.11_2} and \eqref{25.3.11_3} implies 
\begin{equation}\label{lisays}
\ell(g(\gamma_m))^2 \lesssim m^{-q} |g(A_m)| ^{1-\frac{\alpha q}{2}} .
\end{equation}

Let $\sigma = 4/(2-\alpha q ) $ and $\sigma'$ be the conjugate index of $\sigma .$
This  $\sigma$ is well-defined because $q <2 /\alpha .$
Moreover 
$q \sigma' /2 >1$ if and only if $ q >2 /(2-\alpha) .$
By H\"older's inequality, we obtain from \eqref{lisays} that
\begin{align*}
\sum^{\infty} _{m=1} \ell(g(\gamma_m)) \lesssim  & \sum^{\infty} _{m=1} m^{-\frac{q}{2}} |g(A_m)|^{\frac{1}{\sigma}} \\
\le  & \left(\sum^{\infty} _{m=1} m^{-\frac{q \sigma'}{2}} \right)^{\frac{1}{\sigma'}}
\left(\sum^{\infty} _{m=1} |g(A_m)| \right)^{\frac{1}{\sigma}} 
\lesssim | g(\cup^{\infty} _{m=1} A_m) |^{\frac{1}{\sigma}}.
\end{align*}
Set $\Delta_{n,j} = g(\mbd \cap \mbh).$ Then $\Delta_{n,j}$ is a subset of $\Omega,$ and  $\partial \Delta_{n,j}$ is composed of $\Gamma_{n,j}$ and a component from  $\Gamma \setminus \partial \Gamma_{n,j} . $
It follows from $ g(\cup^{\infty} _{m=1} A_m) \subset  \Delta_{n,j}$ that $\ell(\cup^{\infty} _{m=1} g(\gamma_m) ) \lesssim |\Delta_{n,j}|^{1/\sigma } .$ Analogously,  $\ell(\Gamma_{n,j} \setminus \cup^{\infty} _{m=1} g(\gamma_m)) \lesssim |\Delta_{n,j}|^{1/\sigma } .$ Therefore, 
\begin{equation}\label{25.3.12_1}
 \ell(\Gamma_{n,j} )^{\sigma} \lesssim |\Delta_{n,j}| .
\end{equation}
Notice that $f^{-1}(\Gamma_{n,j})$ is the hyperbolic geodesic joining the endpoints of $I_{n,j}.$ Since $f$ extends homeomorphically to $\overline {\mbd}$ by the Caratheodory-Osgood theorem and $\Omega$ has finite area, it follows that $\sup_j |\Delta_{n,j}|\to 0$ when $n\to \infty$ and hence by \eqref{25.3.12_1} this also holds for $\ell(\Gamma_{n,j}).$

Pick $p $ with $ p <\sigma .$
Let $\tau = \sigma /p$ and $\tau'$ be the conjugate exponent of $\tau.$
Then $p-2+1/\tau' <0 $ if and only if $p <4/ (2+\alpha q) .$
Therefore, by H\"older's inequality and \eqref{25.3.12_1}, we have
\begin{align} 
\sum^{\infty}_{n=n_0} 2^{(p-2)n} \sum^{2^n}_{j=1} \ell(\Gamma_{n,j})^p
\le & \sum^{\infty}_{n=n_0} 2^{(p-2)n} \left(\sum^{2^n} _{j=1} 1 \right)^\frac{1}{\tau'} \left(\sum^{2^n} _{j=1} \ell(\Gamma_{n,j})^{\sigma} \right)^{\frac{1}{\tau}}  \notag \\
\lesssim &  \sum^{\infty}_{n=n_0} 2^{(p-2)n} 2^{\frac{n}{\tau '}} |\cup^{2^n} _{j=1} \Delta_{n,j}|^{\frac{1}{\tau}} , \label{26.4.12_1}
\end{align}
where the second inequality comes from $\Delta_{n,j} \cap \Delta_{n,k} =\emptyset $ when $j \neq k .$
Since $\cup^{2^n} _{j=1} \Delta_{n,j} \subset \Omega ,$ we have 
\begin{align}
\sum^{\infty}_{n=n_0} 2^{(p-2)n} 2^{\frac{n}{\tau '}} |\cup^{2^n} _{j=1} \Delta_{n,j}|^{\frac{1}{\tau}} 
\le&  |\Omega|^{\frac{1}{\tau}} \sum^{\infty}_{n=n_0} 2^{n(p-2+\frac{1}{\tau '})} \label{26.5.5_1}\\
=  &  |\Omega| ^{\frac{1}{\tau}} \sum^{\infty}_{n=n_0} 2^{n(p \frac{2+\alpha q}{4} -1)} <\infty \notag
\end{align}
for all $p <4/ (2+\alpha q) ,$ where the equality is from  $\tau = \sigma /p$  and $\sigma = 4/(2-\alpha q )  .$
It follows from \eqref{26.4.12_1} and \eqref{26.5.5_1} that $\sum^{\infty}_{n=n_0} 2^{(p-2)n} \sum^{2^n}_{j=1} \ell(\Gamma_{n,j})^p <\infty$ for all $p <4/ (2+\alpha q) .$ 
Since $ q > 2/(2-\alpha) $  can be taken arbitrarily close to $2/(2-\alpha),$ we conclude that the above holds for any  $p< 2-\alpha .$ The claim follows by applying Theorem \ref{23.12.10_1}.
\end{proof}

The sharpness in Theorem \ref{perus} is a consequence of the following example.

\begin{example}\label{23.12.12_1}
Given $\alpha \in (0,1),$
there is a Jordan domain $\Omega$ with
$$h_{\Omega} (z_1 ,z) \le C d^{-\alpha} (z ,\pao)$$ for every $z \in \Omega ,$ where $z_1 \in \Omega$ is fixed and $C>0$ is a constant, together with
a parameterization $\varphi : \mbs \ra \pao$ failing to have a homeomorphic extension $h: \mbd \ra \Omega$ in the class $W^{1,2-\alpha} (\mbd ,\Omega) .$
\end{example}

Given  $\alpha \in (0,1)$ as in the statement, we fix $\beta \in (1,2)$ and $\tau >0$ with 
\begin{equation}\label{24.4.15_1}
\frac{1-\alpha}{2-\alpha} <\tau <  \frac{1-(\beta -1) \alpha}{2-\alpha}.
\end{equation}
We begin with the construction of  $\Omega .$ By a decomposition of the boundary $\pao ,$ we then obtain a homeomorphic parametrization $\varphi$ of $\pao ,$ showing that $\Omega$ is a Jordan domain.
For simplicity, we define this $\varphi$ on the boundary of a square.
We close by verifying the pointwise quasihyperbolic growth condition and the non-existence of homeomorphic $W^{1,2-\alpha}-$extensions of $\varphi .$

\bigskip
{\it Construction of $\Omega .$}
\bigskip

\begin{figure}[htbp]
\centering

  \includegraphics[width=0.2\textwidth]{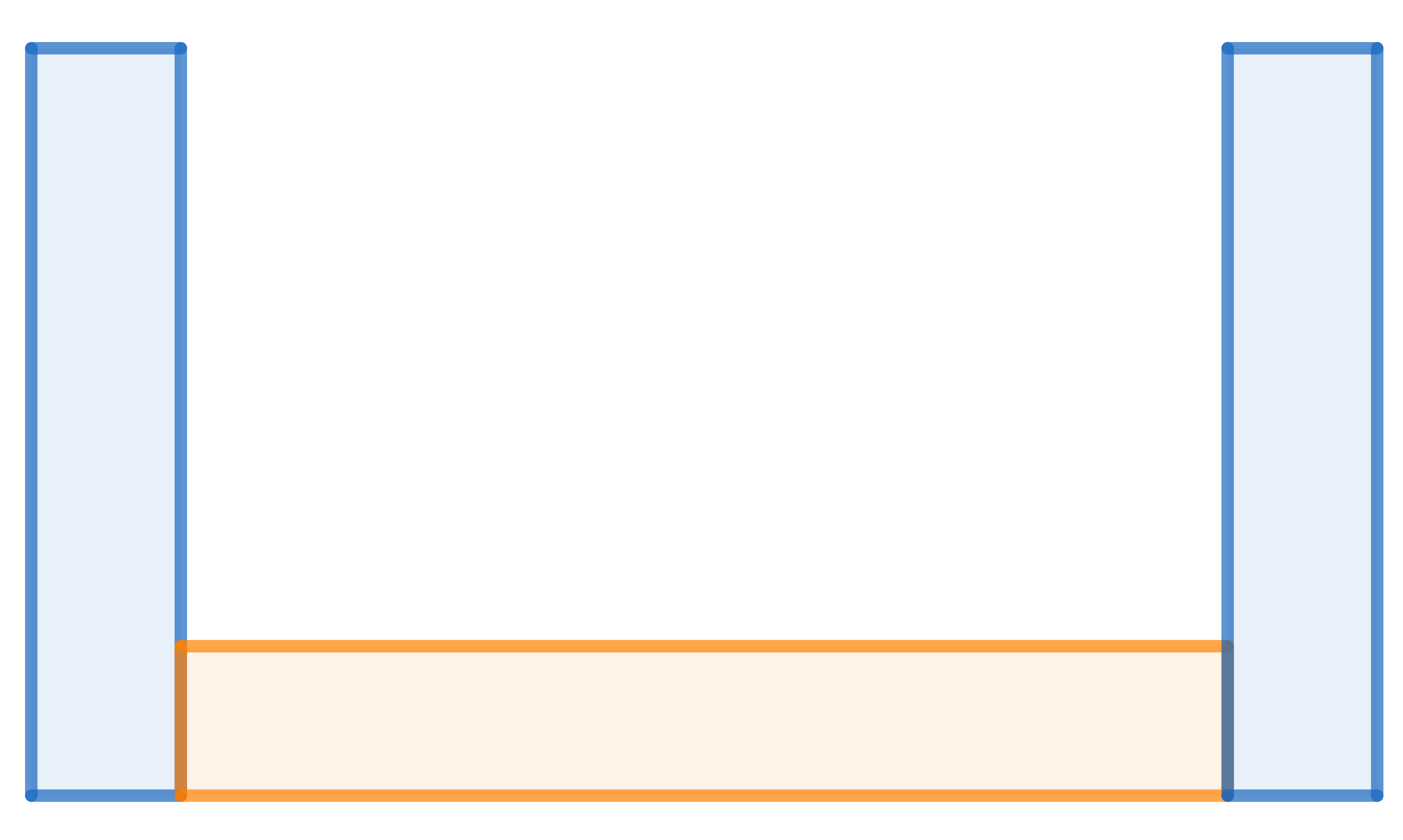}  
  
\caption{The $U$-type tube domain $U_n .$ We require that each arm rectangle contains the square at the corner of $U_n .$}
\label{26.4.9_1}
\end{figure}

\begin{figure}[htbp]
\centering

  \includegraphics{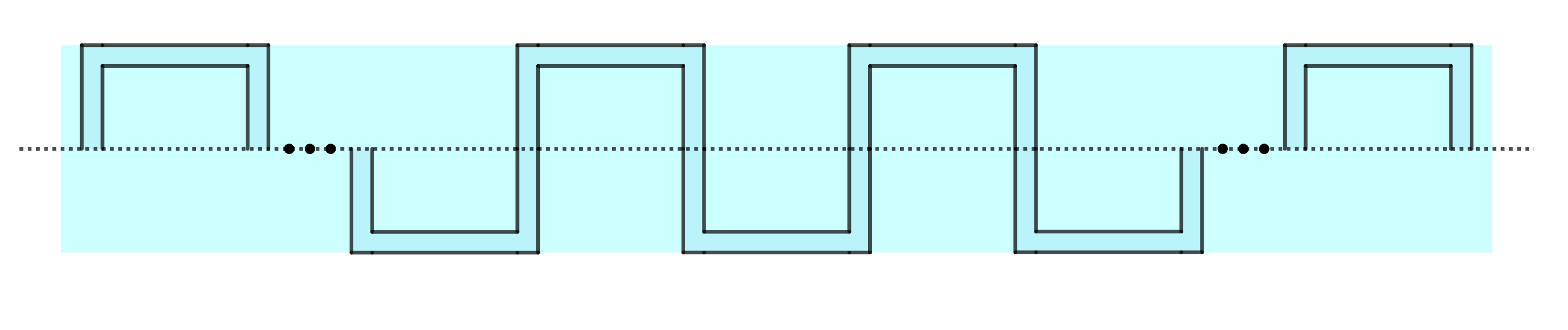}

\caption{ The snake-like domain $S_n$. 
The omitted parts are indicated by dots ($\cdots$). 
$R_n$ is the  coloured part.}
\label{25.2.16_1}
\end{figure}

For $n=1,2,...,$ we let $R_n$ be an open rectangle with side lengths 
\begin{equation}\label{24.5.14_1}
a_n = 2^{-\frac{n}{2}} n^{-\epsilon} \qquad \mbox{ and } \qquad
b_n = 2^{-\frac{n}{2}} n^{-1-2\epsilon} ,
\end{equation}
where $\epsilon $ is positive and satisfies $1+3\epsilon < \tau (2-\alpha)/(1-\alpha) .$
Using $R_n ,$ we construct a snake-like set $S_n .$
See Figure \ref{25.2.16_1} for the set $S_n .$
The length of the middle axis of $S_n$ will be approximately
\begin{equation}\label{24.5.16_7}
\ell_n = 2^{\frac{n(1-\alpha)}{\alpha -2}} n^{-\tau}
\end{equation}
and the diameter of its opening is  $\ell^{1/(1-\alpha)} _n .$ To do this, we first introduce an $U$-type tube domain $U_n .$ See Figure \ref{26.4.9_1}.
It is composed of two vertical (arm) rectangles and a (base) horizontal rectangle. 
These two arm rectangles have side lengths  $b_n /2 $ and
$\ell^{1/(1-\alpha)} _n .$ The base rectangle has side lengths $c_n $ and  $\ell^{1/(1-\alpha)} _n $ with $ c_n \ll a_n $ to be determined.
Flip $U_n$ over to obtain an inverted $U$-type domain, denoted by $U'_n .$
Next, we will glue $U_n $ and $U'_n $ in a suitable manner.
In this process, base segments of $U_n $ and of $U'_n $ are located on the $a_n$-sides of $R_n .$
In the first step, we position $U_n$ so that the $a_n$-sides are perpendicular to the line joining centers of $R_n$ and of the base rectangle of $U_n .$
Secondly,  left from this $U_n ,$ we attach  a copy of  $U' _n$ so that  its right arm is glued to the left arm of $U_n .$ We glue another copy of  $U' _n$ symmetrically to the right arm of $U_n .$ We continue by attaching  another copy of  $U_n$ to the left arm of the leftmost $U'_n$ and another copy of  $U_n$ to the right arm  of the rightmost  $U'_n .$
Inductively, by gluing the total number $\kappa_n$ of copies of $U_n$ and $U' _n$ we obtain a connected snake-type domain $S_n$ with 
$ S_n \subset R_n .$
Of course, the neighbouring $U_n$ and $U' _n$ touch different $a_n$-sides of $R_n .$

The above process shows that $\kappa_n$ is odd. Moreover, we take $\kappa_n$ as the largest odd number with
\begin{numcases}{}
  & $\kappa_n (c_n +b_n ) \le \ell_n ,$ \label{25.2.17_1}\\
& $\kappa_n (c_n + \ell^{\frac{1}{1-\alpha}} _n ) + \ell^{\frac{1}{1-\alpha}} _n \le a_n - \ell^{\frac{1}{1-\alpha}} _n .$ \label{25.2.17_2}
\end{numcases}
The term $c_n +b_n$ in \eqref{25.2.17_1} is the middle axis length of $U_n .$
Hence,  the middle axis length of $S_n$ is approximately $\ell_n .$
The term $\kappa_n (c_n + \ell^{\frac{1}{1-\alpha}} _n ) + \ell^{\frac{1}{1-\alpha}} _n$ in \eqref{25.2.17_2}
is the  length of the projection of $S_n$ to the (lower) $a_n$-side of $R_n .$
Hence, by  \eqref{25.2.17_2} we see that $S_n \subset R_n $ and the distance between $S_n$ and the $b_n$-sides of $R_n$ is at least $\ell^{\frac{1}{1-\alpha}} /2 .$ 
Since $\kappa_n c_n \in (c_n , a_n),$ it follows from \eqref{25.2.17_1} that $\kappa_n \approx \ell_n /b_n .$
By \eqref{25.2.17_2} we have  $\kappa_n (c_n + \ell^{\frac{1}{1-\alpha}} _n ) \approx a_n .$
Together with the assumption $1+3\epsilon < \tau (2-\alpha)/(1-\alpha) ,$ we obtain $c_n \approx a_n /\kappa_n .$

\begin{figure}[htbp]
\centering

  \includegraphics[width=0.99\textwidth]{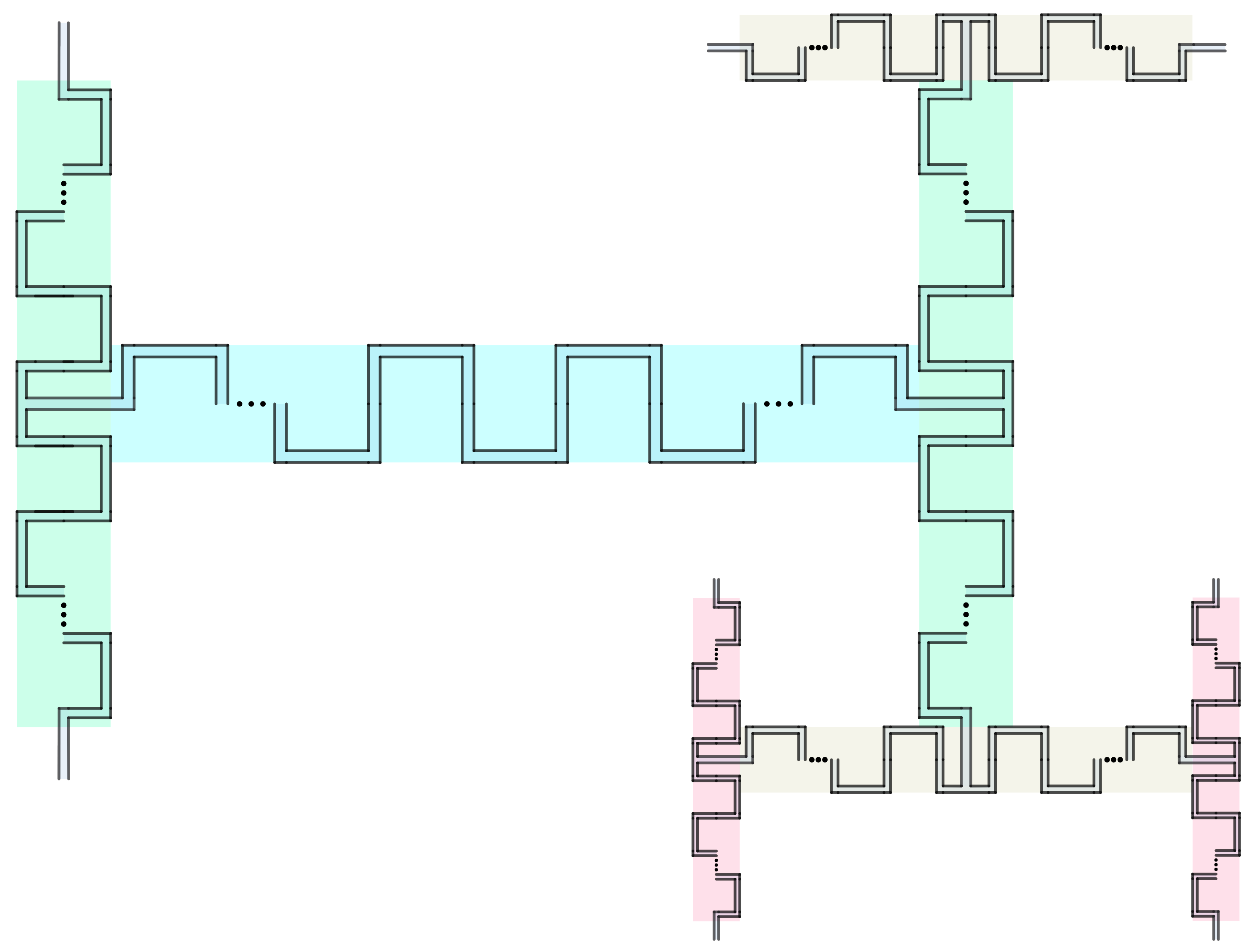}

\caption{For $T_1 ,$ we glue two $T_2 \in \mathfrak{T}_2 $ to the sides of $T_1 .$ For both  $T_2 ,$ we attach two $T_3$ from $\mathfrak{T}_3$ to the sides. In this figure, this is only done for one $T_2 .$ Analogously, we glue two $T_4$ to the  sides of every $T_3 .$}
\label{25.4.8_1}
\end{figure}

We next suitably arrange copies of $S_n$ to obtain a domain $\Omega .$ Since $S_n \subset R_n ,$ it suffices to make a suitable arrangement for the rectangles $R_n .$
Denote by $d_n$ the distance between $S_n$ and the $b_n$-sides of $R_n .$
On both of the  $b_1$-sides of $R_1 ,$
we glue an $R_2$ so that the $a_1$-sides of $R_1 $ are parallel to
the line $L_{1,2}$ joining the centers of $R_1$ and of $R_2 $ and the $a_2$-sides of the $R_2$'s are perpendicular to the $a_1$ sides of $R_1.$
To connect $S_1$ and these copies of $S_2 ,$ we use  rectangular domains $E_1$ whose centers are located in $L_{1,2}$ and whose side lengths are $\ell^{1/(1-\alpha)} _1$ and $\ell^{1/(1-\alpha)} _1 +d_1 +b_2 . $
Then, $E_1 \cap S_1$ is a half of a square with the side length $\ell^{1/(1-\alpha)} _1$ and  $E_1 \cap S_2$ is a rectangle with side lengths $\ell^{1/(1-\alpha)} _1$ and $\ell^{1/(1-\alpha)} _2 .$
Write $$T_1 = S_1 \cup E_1 .$$
In general, after the $n$-step, we attach the number $2^{n}$ of $R_{n+1} $ and obtain the same number of $T_n .$
In the $(n+1)$-step, on each $b_{n+1}$-side of these $R_{n+1}$ we glue an $R_{n+2}$ such that the $a_{n+1}$-sides are parallel to 
the line $L_{n+1, n+2}$ joining the centers of $R_{n+2}$ and of the attached $R_{n+1}.$
To connect $S_{n+1}$ and $S_{n+2} ,$ we use  open rectangles $E_{n+1}$ whose centers are located in $L_{n+1, n+2} ,$ and side lengths are $\ell^{1/(1-\alpha)} _{n+1}$ and $\ell^{1/(1-\alpha)} _{n+1} + d_{n+1} +b_{n+2} . $
Moreover, $E_{n+1} \cap S_{n+1}$ is a half of a square with the side length $\ell^{1/(1-\alpha)} _{n+1}$ and  $E_{n+1} \cap S_{n+2}$ is a rectangle with side lengths $\ell^{1/(1-\alpha)} _{n+1}$ and $\ell^{1/(1-\alpha)} _{n+2} .$
Let $T_{n+1} = S_{n+1} \cup E_{n+1} . $
We continue this procedure ad inductively. Denote by $\mathfrak{T}_n$ the collection of all the above $T_n $ for $n \ge 1 .$ 
Let 
\begin{equation}\label{25.2.17_3}
\Omega = \cup^{\infty} _{n=1} \cup_{T \in \mathfrak{T}_n} T .
\end{equation}
Denote by $\pao$ the boundary of $\Omega .$ See Figure \ref{25.4.8_1} for the first steps of the construction of $\Omega .$

We claim that $\pao$ does not  have self-intersections. 
Since $S_n \subset R_n ,$  it suffices to show that no $R_{n-1} $ intersects an  $R_k$ with $k \ge n .$ By the geometry of $\Omega ,$
this follows if
\begin{equation}\label{25.2.18_1}
a_n  - \sum^{\infty} _{k=1}  (a_{n+2k} - b_{n+2k -1} ) - 2\sum^{\infty} _{k=1}  b_{n+2k +1} > b _{n -1} .
\end{equation}
Here $a_n$ is the larger side length of an $R_n$ attached to $R_{n-1} .$
The term $a_{n+2k} - b _{n+2k -1} $ in the first sum is the $a_{n+2k}$-side length of $R_{n+2k}$ minus the  $b _{n+2k -1}$-side length of $R_{n+2k -1} .$
The second sum is the total length of $b _{n+2k +1}$-side  of $R_{n+2k +1}$ for $k \ge 1 .$ 
Let $D_1 = a_n  - \sum^{\infty} _{k=1} a_{n+2k}$ and $D_2 = 2b_{n+1}- \sum^{\infty} _{k=1}  b_{n+2k -3} .$
The above inequality is equivalent to 
$D_1 +D_2 > 0 .$
By \eqref{24.5.14_1} we estimate 
\begin{align*}
D_1= & 2^{-\frac{n}{2}} \left(n^{-\epsilon} - \sum^{\infty} _{k=1} 2^{-k} (n+2k)^{-\epsilon} \right) 
>  2^{-\frac{n}{2}} \left(n^{-\epsilon} - (n+2)^{-\epsilon} \sum^{\infty} _{k=1} 2^{-k}  \right) \\
= &  2^{-\frac{n}{2}} (n^{-\epsilon} - (n+2)^{-\epsilon} )
\end{align*}
and 
\begin{align*}
D_2 = & 
2b_{n+1} -b_{n -1}  - \sum^{\infty} _{i=1}  b_{n+2i -1} 
>  2b_{n+1} -b_{n -1} -2b_{n+1} = -b_{n -1} .
\end{align*}
It follows that 
\begin{align*}
D_1 +D_2 > & 2^{-\frac{n}{2}}  \left(\frac{1}{n^{\epsilon}}- \frac{1}{(n+2)^{\epsilon}} -\frac{\sqrt{2} }{(n-1)^{1+2\epsilon}} \right) \\
>& 2^{-\frac{n}{2}} \left(\frac{2 \epsilon}{(n+2)^{1+\epsilon}} -\frac{\sqrt{2} }{(n-1)^{1+2\epsilon}} \right).
\end{align*}
For $\epsilon >0 ,$ there is a large enough $n_\epsilon$ such that $D_1 +D_2 >0$ for all $n \ge n_\epsilon .$
Hence, the inequality \eqref{25.2.18_1} holds for all $n \ge n_\epsilon .$  Thus, we start the construction \eqref{25.2.17_3} from $n_{\epsilon} .$
For notational simplicity, we let $n_\epsilon =1 $ in what follows.
This does not have a significant effect on the  argument.

\bigskip
{\it A decomposition of $\pao .$}
\bigskip

We next give a decomposition of $\pao .$
Recall that $\mathfrak{T}_n $ in \eqref{25.2.17_3} is composed of tube-type domains. Every $ T_n$ in $\mathfrak{T}_n $ with $n\ge 2$ meets two elements in $\mathfrak{T}_{n+1} $  and one in $\mathfrak{T}_{n-1} .$  Denote by $A (T_n )$ the union of this $(n-1)$-step set and the closure of these two $(n+1)$-step sets.
Then
\begin{equation*}\label{25.2.20_3}
\gamma(T_n) := \partial T_n \setminus A (T_n )
\end{equation*}
is composed of three arcs: two with the same length and one with a longer length. Moreover, the longest one does not contain its endpoints. Both shorter arcs contain only one endpoint, which is on the boundary of the attached $T_{n-1} .$
For $n=1, $ there is no $T_{n-1}$ and we simply set $\gamma(T_1)$ 
to be the union of the two associated boundary arcs, denoted by
\begin{equation*}\label{25.2.20_1}
\gamma_1 \qquad \mbox{ and } \qquad \gamma_2 .
\end{equation*}
They do not contain their endpoints.
We label them so  that $\gamma_1$ is above  $\gamma_2$ in Figure \ref{25.4.8_1}. Given this $\gamma_1 ,$ we select from $ \cup_{T_{2} \in \mathfrak{T}_{2}} \gamma(T_2)$ two arcs  such that each of them shares one endpoint with $\gamma_1 .$
Denote them by $\gamma^1 _2$ and $\gamma^2 _2 .$
Continue in a similar manner until the step $n.$
Then, we have two arcs $\gamma^1 _n$ and $\gamma^2 _n$ from $\cup_{T_{n} \in \mathfrak{T}_{n}} \gamma(T_n)$ so that $ \gamma_1 \cup \cup^n _{i=1} \gamma^1 _i \cup \gamma^2 _i$ is an open Jordan arc.
Here an open Jordan arc is a homeomorphic image of an open interval.
For this $\gamma^1 _n $ there is an arc $\gamma^1 _{n+1} $ from $\cup_{T_{n+1} \in \mathfrak{T}_{n+1}} \gamma(T_{n+1})$ that shares one endpoint with  $\gamma^1 _n .$ 
Analogously,  there is $\gamma^2 _{n+1} \in \cup_{T_{n+1} \in \mathfrak{T}_{n+1}} \gamma(T_{n+1})$ such that it shares one endpoint with $\gamma^2 _n .$ Continue this process inductively.
Denote 
\begin{equation*}\label{25.2.20_2}
\Gamma_1 = \gamma_1 \cup \cup^{\infty} _{i=2} \gamma^1 _i \cup \gamma^2 _i .
\end{equation*}
Then $\Gamma_1$
is an open Jordan arc in $\pao .$  Moreover, $\Gamma_1$ is symmetric with respect to the perpendicular bisector of $a_1$-sides of $R_1 .$
Replace $\gamma_1$ with $\gamma_2$ and then run again the process for $\Gamma_1 .$ 
This results in  another open Jordan arc $\Gamma_2 \subset \pao .$

Pick  $T_2 \in  \mathfrak{T}_{2} .$ Then $\gamma (T_2) $ is composed of three arcs. The two shorter ones have been chosen as  parts of $\Gamma_1 \cup \Gamma_2 .$ 
We consider the longest one as a new beginning of the operations we performed for $\Gamma_1 .$ In this way, we obtain an open Jordan arc, which is a part of $\pao .$
Since there are two $T_2$ in $\mathfrak{T}_{2} ,$
we have two such open Jordan arcs in total. 
Denote by $\mathfrak{G}' _2$ a singleton, whose only element is the first Jordan arc. Analogously, $\mathfrak{G}'' _2 $ is also  a singleton whose element is the other Jordan arc.
In general, consider $T_n \in  \mathfrak{T}_{n} $ for $n \ge 3 .$
Then $\gamma (T_n)$ is composed of three arcs. The two shorter ones have been applied to construct the previous Jordan arcs.
We again apply the operations used for $\Gamma_1$ to the longest one. This results in a new open Jordan arc. 
This procedure gives   $2^{n-1}$ open Jordan arcs corresponding to $\mathfrak{T}_{n} .$ 
We split these Jordan arcs into two equal groups, which are denoted by \begin{equation*}
\mathfrak{G}' _n \qquad \mbox{ and } \qquad \mathfrak{G}'' _n 
\end{equation*} 
so that each desired $\gamma' \in \mathfrak{G}' _n $ is located on the left half of $T_1 ,$ which is the initial subdomain of the $\Omega.$ 
Similarly, each $ \gamma '' \in \mathfrak{G}'' _n $ is located on the right half of  $T_1 .$
Let $E' $ be the collection of accumulation points of $\cup^{\infty} _{n =2} \cup_{\gamma \in \mathfrak{G}' _n} \gamma $ in  $\cup^{\infty} _{n =2} \cup_{\gamma \in \mathfrak{G}' _n} \gamma .$    We define  a subset $E''$ of  $\cup^{\infty} _{n =2} \cup_{\gamma \in \mathfrak{G}'' _n} \gamma $ analogously. 
Set
\begin{equation}\label{25.2.23_2}
\Gamma_3 = E' \cup  \cup^{\infty} _{n =2} \cup_{\gamma \in \mathfrak{G}' _n} \gamma
\qquad
\mbox{ and }
\qquad
\Gamma_4 = E''   \cup \cup^{\infty} _{n =2} \cup_{\gamma \in \mathfrak{G}'' _n} \gamma .
\end{equation}
This allows us to  decompose $\pao$ as 
\begin{equation*}
\pao = \Gamma_1 \cup  \Gamma_2 \cup  \Gamma_3 \cup  \Gamma_4 .
\end{equation*}

\bigskip
{\it Parametrization $\varphi$ of $\pao .$}
\bigskip

For $n \ge 1,$ let
\begin{equation}\label{24.5.16_1}
\hat{\ell}_n = 2^{-n} n^{-\beta} 
\end{equation}
for a fixed $\beta \in (1,2) .$
Let $Q$ be a square with the side length $\hat{\ell}:=  \sum^{\infty}_{n=1} 2^{n-1} \hat{\ell}_n .$
 We construct a homeomorphic parametrization  $\varphi : \partial Q \ra \pao .$
To begin, we decompose $\partial Q$ as $\partial Q = \cup^4 _{k=1} I_k ,$ where $I_1 = (0, \hat{\ell})\times \{\hat{\ell}\} ,$ 
$I_2 = (0, \hat{\ell})\times \{0\} ,$
$I_3 = \{0\} \times [0, \hat{\ell}]$ 
and $I_4 = \{\hat{\ell}\} \times [0, \hat{\ell}].$  
We first construct $\varphi : I_k \ra \Gamma_k$ for $k=1$ and $2 .$
Then, using decompositions of $I_3$ and $I_4 $  generated by suitable Cantor sets, we construct $\varphi : I_k \ra \Gamma_k$ for $k=3$ and $4 .$
A combination of these maps gives a homeomorphic parametrization of $\pao .$

We introduce a piecewise defined homeomorphism $\varphi$ from $I_1$ onto $\Gamma_1 .$
Recall that 
$\Gamma_1 = \gamma_1 \cup \cup^{\infty} _{i=2} \gamma^1 _i \cup \gamma^2 _i ,$
 where $\gamma_1 ,\ \gamma^1 _i$ and $\gamma^2 _i $ are Jordan arcs.  
Here $\gamma_1$ is an open arc, whose left endpoint is in $ \gamma^1 _2$ and the right one is in $\gamma^2 _2 .$ In general, $\gamma^1 _i$ for $i \ge 3$ shares only one endpoint  with $\gamma^1 _{i-1} .$ The analogous conclusion holds for  $\gamma^2 _i $ and  $\gamma^2 _{i-1} .$
Therefore, $\{ \gamma_1 ,\ \gamma^1 _i ,\ \gamma^2 _i \}$ is a pairwise disjoint collection.
We make a similar decomposition for the open interval $I_1$ as $ I_1 = \lambda_1 \cup \cup^{\infty} _{i=1} \lambda^1 _i \cup \lambda^2 _i .$ The intervals $\{\lambda_1, \ \lambda^1 _i ,\ \lambda^2 _i \}$ have the same properties as the corresponding $\{ \gamma_1 ,\ \gamma^1 _i ,\ \gamma^2 _i \} .$
We define 
\begin{equation}\label{25.2.23_1}
\varphi : I_1 \ra \Gamma_1
\end{equation}
as a piecewise defined homeomorphism.  It maps each component of $I_1$ onto the corresponding one in $\Gamma_1 $ by a homeomorphism of constant speed. 
Analogously, there is a piecewise defined homeomorphism $\varphi : I_2 \ra \Gamma_2 .$

In order to construct a homeomorphism $\varphi : I_3 \ra \Gamma_3 ,$ we introduce a  decomposition for the closed interval $I_3 $ via a Cantor set. Recall that $I_3$ is of length $\sum^{\infty} _{n=1} 2^{n-1} \hat{\ell}_n .$ 
We remove from $I_3 $ a middle open interval $K_1$ of length $\hat{\ell}_1 .$ There are two remaining closed intervals.
We remove middle open intervals of length $\hat{\ell}_2  $ from them. Denote by $K_2$ the collection of these two removed intervals.
In general, after step $n,$ there are  $2^n$ remaining closed intervals. 
From each of them, we remove a middle open interval of length $\hat{\ell}_{n+1}  .$ 
Let $K_{n+1}$ be the collection of all removed intervals in this step.
Then  we are left with $2^{n+1}$ closed intervals with equal length. 
By continuing this process inductively, we obtain a Cantor set which is the intersection of all remaining closed intervals.  Denote it by $\mathcal{C}.$ 
Then
\begin{equation}\label{25.2.26_1}
I_3 = \mathcal{C} \cup \cup^{\infty} _{i=1} \cup_{\lambda \in K_i} \lambda .
\end{equation}
Since the length of $I_3$ is $\sum^{\infty} _{n=1} 2^{n-1} \hat{\ell}_n ,$ we see that $\mathcal{C}$ is of length zero.

Next, we construct a homeomorphism $\varphi$ from $I_3$ onto $\Gamma_3 .$
Recall from \eqref{25.2.23_2} that  $\Gamma_3 = E' \cup  \cup^{\infty} _{n =2} \cup_{\gamma \in \mathfrak{G}' _n}  \gamma .$ 
Analogously to  \eqref{25.2.23_1},  there is a homeomorphism $\varphi$ from $K_1$ onto $\mathfrak{G}' _2 .$ For $i \ge 2,$ $K_i$ consists of  $2^{i-1}$  open intervals  and  $ \mathfrak{G}' _{i+1}$ consists of an equal number of open Jordan arcs. We map each of these intervals onto an open arc in $ \mathfrak{G}' _{i+1}  $ by a similar homeomorphism as in \eqref{25.2.23_1}, respecting the order.
This results in  an order-preserving homeomorphism 
\begin{equation}\label{25.2.23_3}
\varphi (z) : 
\cup_{\lambda \in K_i} \lambda \ra  \cup_{\gamma \in \mathfrak{G}' _{i+1}  } \gamma  
\end{equation}
for $i \ge 2.$
Denote by $\partial K_i$ the collection of the endpoints of the intervals in $K_i .$
The homeomorphism $\varphi$ has a natural homeomorphic extension to $\cup^{\infty} _{i =1} \partial K_i .$
It remains to define $\varphi$ on $\mathcal{C} \setminus \cup^{\infty} _{i =1} \partial K_i .$
For every $z $ in this set, we take $z_i \in \partial K_i$ such that $|z_i -z| = \inf_{\xi \in \cup_{\lambda \in K_i} \lambda} |\xi -z| .$
By our construction, we have
$|\varphi (z_{i+1}) -\varphi(z_i) | \lesssim a_{i} .$
It follows that $\{\varphi (z_i)\}$ is a Cauchy sequence.
So we can define $\varphi(z)$ as 
\begin{equation}\label{25.2.23_4}
\varphi (z) = \lim _{i \ra \infty} \varphi (z_i)  
\end{equation}
for $z \in \mathcal{C} \setminus \cup^{\infty} _{i =1} \partial K_i .$
Moreover, $\varphi (\mathcal{C}) \subset E'$ by the definition of $E' .$
Combining \eqref{25.2.23_3} with \eqref{25.2.23_4} gives a mapping $\varphi$ from $I_3$ into $\Gamma_3 .$

We claim that $\varphi : I_3 \ra \Gamma_3$ is a homeomorphism. It suffices to prove that $\varphi$ is continuous and bijective.

\begin{enumerate}[(1)]
\item

We first prove the continuity of $\varphi$ at every point $z_0 \in I_3 .$
If $z_0 \in \cup^{\infty} _{i=1} K_i ,$ the continuity is trivial by \eqref{25.2.23_3}. 
Next, assume that  $z_0 \in \cup^{\infty} _{i=1} \partial K_i .$ Without loss of generality, we let $z_0 $ be the left endpoint of one interval in $K_{i_0} .$ 
It then suffices to show that $\varphi$ is left-continuous at $z_0 .$
For any $i >i_0 ,$ denote by $z_i$ the closest point in $\partial K_i$ away from $z_0$ with $z_i <z_0 .$
Then $\lim_{i \ra \infty} z_i =z_0$ and $\{ z_i \}$ is increasing with respect to $ i.$ 
By the construction of $\Omega ,$  the line segment joining $\varphi (z_0)$ to  $\varphi (z_i)$ is covered by a square with the side length $\sum^{\infty}_{n=i-1} a_n .$  Moreover, this square also covers the line segment joining $\varphi (z_0)$ to $\varphi (z)$ for any $z \in [z_i ,z_{i+1}] .$ It follows that 
\begin{equation}\label{25.2.26_2}
|\varphi(z_0) - \varphi(z) | \le \sum^{\infty}_{n=i-1} a_n
\end{equation}
for any $z \in [z_i ,z_{i+1}] .$ Note that $ \lim_{i \ra \infty}  \sum^{\infty}_{n=i-1} a_n=0 .$
It follows from \eqref{25.2.26_2} that $\varphi$ is left-continuous at $z_ 0 .$
The remaining case is that $z_0 \in \mathcal{C} \setminus \cup^{\infty} _{i =1} \partial K_i .$ We can apply the above arguments to prove that $\varphi$ is both left-continuous and right-continuous at $z_0 .$
In conclusion,  $\varphi$ is continuous at every point $z_0 \in I_2 .$

\item
To show that $\varphi $ is injective on $I_3 ,$ we take $\hat{z} ,\  \tilde{z} \in I_3$ with $\hat{z}  \neq \tilde{z}  .$ It follows from \eqref{25.2.18_1} that Jordan arcs in $\cup^{\infty}_{n=2}\mathfrak{G}' _{n}$ are piecewise disjoint.
Hence, if  $\hat{z} ,\  \tilde{z} \in \cup^{\infty} _{i=1} \cup_{\lambda \in K_i} \lambda $ we obtain from \eqref{25.2.23_3} that $\varphi (\hat{z}) \neq \varphi (\tilde{z}) .$
Next, assume that $\hat{z} \in \mathcal{C}$ and $\tilde{z}$ is  in  $K_{i_0} .$
By \eqref{25.2.23_4} and our construction,, there exist $\hat{z}_i \in \cup_{\lambda \in K_i} \lambda $ such that $\varphi (\hat{z}_i)$ belongs to  one of the  $R_{i+1}$ and $\varphi (\hat{z}) = \lim_{i \ra \infty}\varphi (\hat{z}_i) .$ 
Moreover, by \eqref{25.2.23_3} there is $n_0 \ge i_0 +1$ and  a rectangle $R_{n_0}$ such that  $\varphi(\tilde{z}) \in \overline{R_{n_0}} .$
Denote by $e_{n-1}$ the left-hand side of \eqref{25.2.18_1}.
For this $R_{n_0} ,$ by \eqref{25.2.18_1} we see that the distance between $R_{n_0}$ and $R_i$ for any $i \ge n$  is larger than $c_{n_0} ,$ where $c_{n_0}= \min \{(e_{n_0} -b_{n_0} )/2 , (e_{n_0 +1} -b_{n_0 +1} )/2 \} .$
Here we ignore the two $R_{n_0 +1}$ attached to this $R_{n_0} .$
Then
\begin{equation*}
|\varphi (\hat{z}) - \varphi(\tilde{z})| 
= \lim_{i \ra \infty} |\varphi (\hat{z}_i) - \varphi(\tilde{z})| 
\ge c_{n_0} .
\end{equation*}
The remaining case is $\hat{z} ,\  \tilde{z}  \in \mathcal{C} .$
Assume that $ \hat{z} < \tilde{z} ,$ where we identify $I_3$ with an interval on the real axis.
For any $\lambda \in \cup^{\infty} _{i=1} K_i ,$ there is either $\lambda \cap ( \hat{z} , \tilde{z} ) = \emptyset$ or $\lambda \cap ( \hat{z} , \tilde{z} ) = \lambda .$
Pick  $\lambda_m $ from $K_m$ with the minimal $m$ such that $\lambda_m \subset ( \hat{z} , \tilde{z} ) .$ There is only one such $\lambda_m$ in this $K_m .$
From the left of $\lambda_m, $ we pick the closest  $\lambda_{m_1} \in \{\lambda \in K_i : 1 \le i \le m-1\}$
to $\lambda_m .$  Analogously, we pick  $\lambda_{m_2} \in \{\lambda \in K_i : 1 \le i \le m-1\} $ from the right of $\lambda_m .$ Denote by $z_1$ the right endpoint of $\lambda_{m_1}$ and by $z_2$ the left endpoint of $\lambda_{m_2} .$ Then 
\begin{equation*}
\lambda _m \subset [\hat{z} , \tilde{z}] \subset [z_1 ,z_2].
\end{equation*}
and $\lambda_m$ is a middle interval of $[z_1 ,z_2]. $ 
Let  $x ,y \in ([z_1, z_2] \setminus \lambda_m )\cap (\cup^{\infty} _{i=m+1} \cup_{\lambda \in K_i} \lambda)$ with $x < y .$
By \eqref{25.2.18_1} and the geometry of our construction, we obtain that $|\varphi (x) -\varphi (y)| > b_{m-1} .$
This combined with  \eqref{25.2.23_4} gives $|\varphi (\hat{z}) - \varphi (\tilde{z} )| >0 .$
Therefore, $\varphi$ is injective.

\item
We show that $\varphi : I_3 \ra \Gamma_3$ is surjective.
By  \eqref{25.2.23_2} and  \eqref{25.2.23_3}, it suffices to show that each $w \in E'$ has a preimage in $I_3 .$ For this $w ,$ there is $\{w_n\} \subset \cup^{\infty} _{n=2}  \cup_{\gamma \in \mathfrak{G}' _n} \gamma$ such that 
$\lim_{n \ra \infty} w_n = w. $
By \eqref{25.2.23_3}, we obtain $\{z_n\} \subset \cup^{\infty} _{i=1} \cup_{\lambda \in K_i} \lambda $  such that $\varphi (z_{n}) = w_n .$  
There is a subsequence $\{\hat{z}_{i} \}$ of $\{z_{n} \}$ and  $z \in I_3$ such that $\lim_{i \ra \infty} \hat{z}_{i} =z .$ Since $\varphi$ is continuous on $I_3 ,$ we obtain that 
\begin{equation*}
w = \lim_{i \ra \infty} \varphi (\hat{z}_{i} ) = \varphi (z).
\end{equation*}
Therefore, for every $ w \in E'$ there is $z \in I_3$ such that $\varphi (z) =w .$

\end{enumerate}

Analogously to $\varphi :  I_3 \ra \Gamma_3 ,$
 there is a homeomorphism from $I_4$ onto $\Gamma_4 .$

\bigskip
{\it The estimate $h_{\Omega} (z_1 ,z) \le C d^{-\alpha} (z ,\pao) .$} 
\bigskip

Every $S_i$ for $i \ge 1$ is composed of $U_i$ and its variants. For the original $U_i$ in $S_i ,$ we denote by $z_i$ the center of its base rectangle.
Fix $z \in \Omega.$ By \eqref{25.2.17_3}, we know that $z \in T_n $ for some $n \ge 1 .$
Then $d (z , \partial T_n) = d(z ,\pao) .$
Let $\hat{z}$ be the projection of $z$ on the middle axis of $T_n .$
Denote by $\lambda (z_n, \hat{z})$ the  middle axis in $T_n$ joining $z_n$ to $\hat{z} .$
Then $\ell(\lambda (z_n, \hat{z})) \lesssim \ell_n ,$ where $\ell_n$ is from \eqref{24.5.16_7}.
We have
\begin{equation}\label{25.3.3_1}
h_{\Omega} (z_n, \hat{z}) 
\le  \int_{\lambda (z_n, \hat{z})} \frac{d \xi}{d(\xi , \partial T_n)}
\lesssim  
\ell^{\frac{-\alpha}{1-\alpha}} _n 
 \le d^{-\alpha} (z ,\pao ) .
\end{equation}
Let $\gamma(\hat{z} ,z) \subset T_n$ be the line segment joining $\hat{z}$ to $z .$
Then 
\begin{equation*}
h_{\Omega} (\hat{z}, z) \le  \int_{\gamma(\hat{z} ,z)  } \frac{d \xi}{d(\xi ,\partial T_n)} \le  \ln \frac{\ell^{\frac{1}{1-\alpha}} _n }{d (z ,\partial T_n )} \lesssim  d^{-\alpha} (z ,\pao ) .
\end{equation*}
Hence,
$h_{\Omega } (z_n, z) \le   h_{\Omega  } (z_n, \hat{z}) + h_{\Omega  } (\hat{z}, z) \lesssim d^{-\alpha} (z ,\pao ) .$
Analogously to \eqref{25.3.3_1}, we have $h_{\Omega} (z_k ,z_{k+1}) \lesssim \ell^{-\alpha /(1-\alpha)} _k $ for all $k \ge 1 .$
Therefore,
\begin{align*}
h_{\Omega} (z_1, z) 
\le & \sum^{n-1} _{k=1} h_{\Omega} (z_k, z_{k+1}) +h_{\Omega } (z_n ,z) 
\lesssim  \sum^{n-1} _{k=1} \ell^{\frac{-\alpha}{1-\alpha}} _k + d^{-\alpha} (z ,\pao  ) \\
\lesssim &  \ell^{\frac{-\alpha}{1-\alpha}} _n +  d^{-\alpha} (z ,\pao ) \\
\lesssim & d^{-\alpha} (z ,\pao) .
\end{align*}
The desired inequality is obtained.

\bigskip
{\it  The failure of the existence of a homeomorphic extension in $W^{1,2-\alpha} (Q,\Omega) .$}
\bigskip

Suppose that  $h : Q \ra \Omega$ is a homeomorphic Sobolev extension of $\varphi .$
We apply  \eqref{25.2.26_1} towards a lower bound for  the integral of $|Dh|^p,$ $p>1,$ over a part of $Q$ close to  $I_3 .$ 
Fix an open interval $\lambda _n \in K_n$ and denote by $p_n$ the middle point of $\lambda_n .$ Let  $\lambda_{n-1} \in K_{n-1}$ be the one in $K_{n-1}$ nearest to $\lambda_n $ and  choose $\lambda_{n+1}$ in $K_{n+1} $  analogously.
Let $r_n$ be the distance between $p_n$ and the middle point of $\lambda_{n+1} ,$ and let $\hat{r}_n$ be the distance between $p_n$ and $\lambda_{n-1} .$ Set $A_n$ to be the intersection of $Q$ and the open annulus centered at $p_n $ with radii $r_n$ and $\hat{r}_n .$ 
By \eqref{24.5.16_1}, we have
\begin{equation}\label{24.5.16_6}
r_n = \frac{\hat{\ell}_n}{2} + \frac{\hat{\ell}_{n+1}}{2} + \sum^{\infty} _{k=n+2} \hat{\ell}_k 2^{k-n-2}  \approx 2^{-n} n^{1-\beta}
\end{equation}
and 
\begin{equation}\label{24.5.16_3}
\hat{r}_n =\frac{\hat{\ell}_n}{2}  + \hat{\ell}_{n+1}  +2\sum^{\infty} _{k=n+2} \hat{\ell}_k 2^{k-n-2}  \approx 2 r_n .
\end{equation}
For $r \in (r_n ,\hat{r}_n),$ we denote by $\zeta_r $ the half circle in $A_n$ with  center $p_n$ and radius $r .$
Since $h$ is a homeomorphic Sobolev mapping,  $h(r \cos  \theta , r \sin \theta)$ is absolutely continuous with respect to $\theta \in [0,\pi]$ for a.e. $r \in (r_n ,\hat{r}_n) .$ Hence
\begin{equation}\label{26.4.10_2}
\ell(h(\zeta_r)) = \int^{\pi} _{0} \left | \frac{\partial}{\partial \theta} h(r \cos  \theta , r \sin \theta) \right| \, d \theta
\le \int_{\zeta_r} |Dh| ds ,
\end{equation}
where the inequality is from a change of variable, and $ds$ is the arc length element over $\zeta_r . $
 Moreover, by the geometry of $\Omega ,$ the curve $h(\zeta_r) $ has to pass through some $T_n \in \mathfrak{T}_n .$ Recall that the length of the middle axis of $T_n$ is comparable to $\ell_n .$ 
So $\ell_n \lesssim \ell(h(\zeta_r)) .$
Together with Jensen's inequality, for $p \in (1,2-\alpha)$ we obtain from \eqref{26.4.10_2} that
\begin{equation*}\label{24.5.16_2}
\ell^p _n r^{1-p} \lesssim \ell(h(\zeta_r))^p r^{1-p}   \le \int_{\zeta_r} |Dh|^p ds .
\end{equation*}
By integrating both sides with respect to $r \in (r_n ,\hat{r}_n),$ we obtain
\begin{equation}\label{24.5.16_4}
\ell^p _n r^{2-p} _n  \approx \int^{\hat{r}_n} _{r_n} \ell^p _n r^{1-p} dr
\lesssim  \int^{\hat{r}_n} _{r_n} \int_{\zeta_r} |Dh|^p ds  dr = \int_{A_n} |Dh|^p dxdy ,
\end{equation}
where  the equality comes from Fubini's theorem.

Let $\sigma_n $ be the maximal square embedded in $A_n$ such that $A_n \setminus \sigma_n$ is symmetric.
Then $p_n$ is the projection of the center of $\sigma_n$ to $\lambda_n .$
There is an absolute constant $c>1$  such that $A_n \subset c \sigma_n ,$ where $c\sigma_n$ the  square concentric to $\sigma_n$ and of diameter $c$ times that of $\sigma_n .$
It follows from \eqref{24.5.16_4} that 
\begin{equation}\label{24.5.16_5}
\ell^p _n r^{-p} _n \lesssim \mint{-}_{A_n} |Dh|^p \lesssim \mint{-}_{c\sigma_n} |Dh|^p \le M_{|Dh|^p } (z) 
\end{equation}
for all $z \in \sigma_n ,$ where $M_{|Dh|^p }$ is the maximal function of the zero extension of  $|Dh|^p $ from $Q$ to the entire plane.
We raise   \eqref{24.5.16_5} to the power $(2-\alpha)/p$ and integrate over $z \in \sigma_n$ to obtain 
\begin{equation}\label{25.3.4_1}
\ell^{2-\alpha} _n r^{\alpha} _n \approx \int_{\sigma_n} (\ell^p _n r^{-p} _n )^{\frac{2-\alpha}{p}} dz \lesssim \int_{\sigma_n} M_{|Dh|^p } (z)^{\frac{2-\alpha}{p}} dz .
\end{equation}
Let $\Sigma_n$ be the collection of all such squares $\sigma_n$ corresponding to $\lambda_n \in K_n .$ 
For any $\sigma_i \in \Sigma_i$ and $\hat{\sigma}_j \in \Sigma_j ,$
we claim that $\sigma_i \cap  \hat{\sigma}_j = \emptyset .$
In fact, if $i=j $ we can prove it by the definition of $\hat{r}_i .$  If $i \neq j ,$ it suffices to consider the case $j =i+1 .$ Then $ \sigma_i \cap \hat{\sigma}_{i+1} = \emptyset$ by $\hat{r}_{i+1} <r_i .$ Since $\Sigma_n$ consists of  $2^{n-1}$ elements, we obtain from \eqref{25.3.4_1} that
\begin{align}
\sum^{\infty} _{n=1} 2^{n-1} \ell^{2-\alpha} _n r^{\alpha} _n 
 \lesssim & \sum^{\infty} _{n=1} \sum_{\sigma \in \Sigma_n}\int_{\sigma} M_{|Dh|^p } (z)^{\frac{2-\alpha}{p}} dz  \notag\\
 \le & \int_{\mbr^2} M_{|Dh|^p } (z)^{\frac{2-\alpha}{p}} dz 
\lesssim  \int_{Q} |Dh|^{2-\alpha}  (z) dz, \label{25.3.4_2}
\end{align}
where the final inequality is from the $L^{(2-\alpha)/p} (\mbr^2)$-boundedness of the maximal operator.
By inserting the values of  $\ell_n$ and $r_n $ from \eqref{24.5.16_7} and \eqref{24.5.16_6} we conclude via  \eqref{24.4.15_1} that
\begin{equation*}
\sum^{\infty} _{n=1} 2^{n-1} \ell^{2-\alpha} _n r^{\alpha} _n 
=\frac{1}{2} \sum^{\infty} _{n=1} n^{\tau (\alpha -2) +(1-\beta) \alpha}
=\infty.
\end{equation*}
Hence, it follows from \eqref{25.3.4_2} that $h$ is not in $W^{1,2-\alpha} (Q,\Omega) .$



\begin{proof}[Proof of Theorem \ref{perus2}]
We fix a dyadic $I$ as in the beginning of the proof of the positive part of Theorem \ref{perus} and notice that the construction of $\{\Gamma_{n,j}\}$ does not depend on $\alpha .$ Hence \eqref{26.5.5_1} allows us to apply Theorem \ref{23.12.10_1} for all $p<2 .$
\end{proof}

\section{Proof of Theorem \ref{26.7.20}}\label{26.8.11_1}

We begin with a lemma.

\begin{lemma}\label{26.7.24_3}
Let $\Omega$ be a Jordan domain,  $\alpha \in (0,1)$ and $\beta >2 .$ Suppose that there is a constant $C$ such that the hyperbolic distance $\rho_{\Omega}$ satisfies the growth condition
\begin{equation}\label{26.7.20_1}
\rho_{\Omega}(z,z_0)\le C d(z,\partial  \Omega)^{-\alpha} \log^{-\beta} \left(\frac{1}{d(z,\partial  \Omega)} \right)
\end{equation}
for some fixed $z_0\in \Omega$ and all $z\in \Omega .$
There exists a constant $M= M(\Omega)$ such that every homeomorphism $\varphi : \mbs \ra \pao $ admits a homeomorphic extension $h:\overline{\mbd}\to \overline{\Omega}$ that belongs to
$W^{1,2-\alpha}(\mbd,\Omega)$ with norm bounded by $M .$
\end{lemma}

\begin{proof}
The proof is analogous to the positive part of Theorem \ref{perus}.  So we follow the notation and main ideas introduced earlier.
Fix a dyadic family $I = \{I_{n,j} \},$ $n=1,2,\dots,$ $j=1,\dots, 2^n,$ of closed arcs on $\partial \mbd .$
Denote by $\Gamma_{n,j}$ the hyperbolic geodesic connecting the endpoints of $\varphi (I_{n,j}) .$
We take the conformal mapping  $g : \mbh \ra \Omega $ from the beginning of the proof of the positive part of Theorem \ref{perus}.
Then $\gamma_m$ is the preimage of sub-arc of  $\Gamma_{n,j}$ under $g .$
Moreover, $A_m$ is a Whitney square associated with $\gamma_m .$

Denote $\sigma = (2-\alpha) /(1-\alpha).$ Let $\sigma' =2-\alpha$ be the conjugate index of $\sigma .$
By imitating  the proof of \eqref{25.3.11_3} for $q= 2/\sigma' ,$
we derive from the assumption \eqref{26.7.20_1} that 
\begin{align*}
\int_{g(A_m)} h^{q} _{\Omega} (z_0 ,\xi) d \xi
\lesssim & \int_{g(A_m)} d^{-\alpha q} (\xi ,\pao) \log^{-\beta q} \left( \frac{1}{d(\xi ,\pao) }\right) d \xi \\
\approx & |g(A_m)|^{\frac{2-\alpha q}{2}} \log^{-\beta q} \left( \frac{1}{|g(A_m)|}\right) .
\end{align*}
Here the quasihyperbolic metric $ h _{\Omega} (z_0 ,\cdot)$ is comparable to the hyperbolic metric $ \rho _{\Omega} (z_0 ,\cdot) .$
Together with \eqref{25.3.11_2}, we then arrive at 
\begin{equation}\label{26.7.21_3}
\ell(g(\gamma_m)) \lesssim m^{- \frac{1}{\sigma'}} a_m ,
\end{equation}
where $a_m = |g(A_m)|^{1/\sigma} \log^{-\beta/\sigma'} ( 1/|g(A_m)| ) . $
Take a pair of complementary $N-$functions $\Phi_1$ and $\Psi_1$ such that 
\begin{equation*}
\Phi_1 (t) \approx t^{\sigma'} \left(\log \frac{1}{t} \right)^{\eta} \qquad \mbox{ and } \qquad 
\Psi_1 (s) \approx s^{\sigma} \left(\log \frac{1}{s} \right)^{\frac{\eta}{1-\sigma'} } .
\end{equation*}
whenever $t \in [0,1)$ and $s \in [0,1) ,$ where $\eta<-1 $ satisfies $\beta +\eta >1 .$ 
This $\eta$ can be selected because $\beta >2 .$
It follows from \eqref{26.7.21_3}  that 
\begin{equation}\label{26.7.21_2}
\sum^{\infty} _{m=1} \ell(g(\gamma_m)) \lesssim  \sum^{\infty} _{m=1} m^{- \frac{1}{\sigma'}} a_m 
\lesssim  ||\{m^{- \frac{1}{\sigma'}} \}||_{\Phi_1} \, ||\{ a_m \} ||_{\Psi_1} ,
\end{equation}
where the inequality at the right-hand side is derived by applying  the generalized H\"older's inequality \eqref{26.7.21_4} to
$\Phi_1$ and $\Psi_1$.
We estimate $||\{ a_m \} ||_{\Psi_1} $ first, and then estimate $ ||\{m^{- 1 /\sigma'} \}||_{\Phi_1}  .$
Let $\lambda= \left( \sum^{\infty}_{k=1} |g(A_k)| \log^{(\beta +\eta)/(\alpha -1)} (1/ \sum^{\infty}_{k=1} |g(A_k)| ) \right)^{1/\sigma}.$ 
Then 
\begin{align}\label{26.7.21_5}
\frac{1}{\lambda^{\sigma}} = & \frac{1}{ \sum^{\infty}_{k=1} |g(A_k)|} \log^{\frac{\beta +\eta}{1-\alpha} } \left(\frac{1}{\sum^{\infty}_{k=1} |g(A_k)|} \right) \notag \\
\le & \frac{1}{ \sum^{\infty}_{k=1} |g(A_k)|} \log^{\frac{\beta +\eta}{1-\alpha} } \left(\frac{1}{|g(A_m)|} \right) 
\end{align}
for any $m \ge 1 .$
Notice that $\lambda \ll 1$ when the generation $n$ of the hyperbolic geodesic $\Gamma_{n,j}$ is sufficiently large.
Therefore, by the formula of $a_m$ we have 
\begin{equation}\label{26.7.21_6}
\left(\log \frac{ \lambda}{a_m} \right)^{\frac{\eta}{1-\sigma'} } \le \left(\log \frac{1}{a_m} \right)^{\frac{\eta}{1-\sigma'} } \approx \left(\log \frac{1}{|g(A_m)|} \right)^{\frac{\eta}{1-\sigma'} }  .
\end{equation}
We derive from the formula of $\Psi_1 $ that
\begin{align}
\sum^{\infty}_{m=1}  \Psi_1 \left(\frac{a_m}{\lambda}  \right) 
\approx & \frac{1}{\lambda^{\sigma}} \sum^{\infty}_{m=1} a^{\sigma} _m \left(\log \frac{ \lambda}{a_m} \right)^{\frac{\eta}{1-\sigma'} } \notag \\
\lesssim & \frac{1}{\sum^{\infty}_{k=1} |g(A_k)| }  \sum^{\infty}_{m=1}  a^{\sigma}_m \log^{\frac{\beta +\eta}{1-\alpha} } \left(\frac{1}{|g(A_m)|} \right)  \log^{\frac{\eta}{1-\sigma'}} \left(\frac{1}{|g(A_m)|} \right) \notag \\
= & \frac{1}{\sum^{\infty}_{k=1} |g(A_k)| }  
\sum^{\infty}_{m=1} |g(A_m)|  =1 ,   \label{26.7.22_1}
\end{align}
where the first inequality is from  \eqref{26.7.21_5} and \eqref{26.7.21_6}.
To obtain the forthcoming equality,  we replace $a _m$ with its given value.
By substituting the values of exponents $\sigma$ and $\sigma'  ,$
we see that the sum of exponents with respect to $\log$ terms is zero. So we obtain the equality on the left-hand side of \eqref{26.7.22_1}.
By the definition of the Luxemburg norm, we then derive from \eqref{26.7.22_1} that $||\{ a_m \} ||_{\Psi_1} \lesssim \lambda .$
To estimate $||\{m^{- 1 /\sigma' } \}||_{\Phi_1} ,$ we notice that 
\begin{equation*}
\sum^{\infty} _{m=1} \Phi_1 (m^{- \frac{1}{\sigma'}}) \approx \sum^{\infty} _{m=1} m^{-1} (\log m)^{\eta} <\infty ,
\end{equation*}
where $\eta <-1 .$
It follows that $||\{m^{- 1 /\sigma' } \}||_{\Phi_1}  \lesssim 1 .$
Therefore,  we derive from \eqref{26.7.21_2} that 
\begin{equation}\label{26.7.22_2}
\ell (\cup^{\infty} _{m=1} g(\gamma_m) ) = \sum^{\infty} _{m=1} \ell(g(\gamma_m)) \lesssim \lambda ,
\end{equation}
where the equality holds because $\{g(\gamma_m)\}$ are piecewise disjoint.
As in the proof of the positive part of Theorem \ref{perus}, we set $\Delta_{n,j} = g(\mbd \cap \mbh).$ Then $\Delta_{n,j}$ be a subset of $\Omega,$ and  $\partial \Delta_{n,j}$ is composed of $\Gamma_{n,j}$ and a component from  $\Gamma \setminus \partial \Gamma_{n,j} . $
Moreover, $ g(\cup^{\infty} _{m=1} A_m) \subset  \Delta_{n,j} .$ 
By replacing $\lambda$ with its given value,
we obtain from \eqref{26.7.22_2} that 
\begin{align}
\ell (\cup^{\infty} _{m=1} g(\gamma_m) ) \lesssim &   \left( \sum^{\infty}_{k=1} |g(A_k)| \log^{\frac{\beta +\eta}{\alpha -1}} \left(\frac{1}{\sum^{\infty}_{k=1} |g(A_k)|} \right) \right)^{\frac{1}{\sigma}} \notag\\
= & |g(\cup^{\infty} _{m=1} A_m)|^{\frac{1}{\sigma}} \log^{\frac{\beta +\eta}{\alpha -2}} \left(\frac{1}{|g(\cup^{\infty} _{m=1} A_m)| } \right) \notag \\
\le &|\Delta_{n,j}|^{\frac{1}{\sigma}} \log^{\frac{\beta +\eta}{\alpha -2}} \left(\frac{1}{|\Delta_{n,j}| } \right) ,\label{26.7.22_3}
\end{align}
where the equality holds because $\{g(A_m)\}$ are piecewise disjoint. 
Let $A(t) = t^{1/\sigma} \log^{(\beta +\eta)/(\alpha -2)} (1/t).$
Since $A (t)$ is strictly increasing,  from $ g(\cup^{\infty} _{m=1} A_m) \subset  \Delta_{n,j} $  we then obtain the inequality in \eqref{26.7.22_3}.
The upper bound in \eqref{26.7.22_3} can also be applied to estimate the length of $\Gamma_{n,j} \setminus \cup^{\infty} _{m=1} g(\gamma_m) .$
Therefore, 
\begin{equation}\label{26.8.11_2}
\ell(\Gamma_{n,j})\lesssim A(|\Delta_{n,j}|).
\end{equation}
Notice that $f^{-1}(\Gamma_{n,j})$ is the hyperbolic geodesic joining the endpoints of $I_{n,j},$ where $f : \mbd \ra \Omega $ is the conformal mapping  as in the proof of the positive part of Theorem \ref{perus}.
 Since $f$ extends homeomorphically to $\overline {\mbd}$ by the Caratheodory-Osgood theorem and $\Omega$ has finite area, it follows that $\sup_j |\Delta_{n,j}|\to 0$ when $n\to \infty$ and hence by \eqref{26.8.11_2} this also holds for $\ell(\Gamma_{n,j}).$ 
Therefore, the additional assumption in Theorem \ref{23.12.10_1} is satisfied. 
Notice that the inverse function $A^{-1} (s)$ has the expression $A^{-1} (s)\approx s^{\sigma} \log^{(\beta +\eta)/(1-\alpha)}(1/s) $ as $s \ll 1 .$
It follows from \eqref{26.8.11_2} that 
\begin{align}
\ell(\Gamma_{n,j})^{\sigma} \log^{\frac{\beta +\eta}{1-\alpha}} \left(\frac{1}{\ell(\Gamma_{n,j})} \right) 
\approx & A^{-1} (\ell(\Gamma_{n,j})) \notag \\
\lesssim & A^{-1} (A(|\Delta_{n,j}|))
= |\Delta_{n,j}| .   \label{26.7.22_8}
\end{align}

Denote $\kappa = -(\beta +\eta)/\alpha ,$ where $\eta$ is as above.
Take a pair of complementary $N-$functions $\Phi_2$ and $\Psi_2$ such that 
\begin{equation*}
\Phi_2 (t) \approx t^{\frac{1}{\alpha}} \left(\log \frac{1}{t} \right)^{\kappa} \qquad \mbox{ and } \qquad 
\Psi_2 (s) \approx s^{\frac{1}{1-\alpha}} \left(\log \frac{1}{s} \right)^{\frac{\kappa \alpha}{\alpha -1} } 
\end{equation*}
whenever $t \in [0,1)$ and $s \in [0,1) .$
By the generalized H\"older's inequality as in \eqref{26.7.20_3},
we have
\begin{equation}\label{26.7.22_10}
\sum^{2^n}_{j=1}\ell(\Gamma_{n,j})^{2-\alpha} \lesssim  ||\{1\}||_{\Phi_2} \, ||\{\ell(\Gamma_{n,j})^{2-\alpha}\}||_{\Psi_2}  .
\end{equation}
We proceed to estimate $||\{1\}||_{\Phi_2}$ and $||\{\ell(\Gamma_{n,j})^{2-\alpha}\}||_{\Psi_2}$ respectively.
Let $\lambda_1 = 2^{\alpha n} \log^{\alpha \kappa } (2^n) .$ By the formula of $\Phi_2 ,$ we then have 
\begin{align*}
\sum^{2^n}_{j=1} \Phi_2(\frac{1}{\lambda_1})\approx & 
\lambda^{-\frac{1}{\alpha}} _1 \sum^{2^n}_{j=1} \log^{\kappa} (\lambda_1) \\
= & 2^{-n} \log^{-\kappa} (2^n) \sum^{2^n}_{j=1} \log^{\kappa} (2^{\alpha n} \log^{ \alpha \kappa} (2^n) ) \lesssim 1 .
\end{align*}
It follows from the definition of the Luxemburg norm that 
\begin{equation}\label{26.7.22_7}
||\{1\}||_{\Phi_2} \lesssim \lambda_1 = 2^{\alpha n} \log^{ \alpha \kappa } (2^n) .
\end{equation}
Let $\lambda_2 = (\sum^{2^n} _{j=1} |\Delta_{n,j}|)^{1-\alpha} .$
It follows from the formula of $\Psi_2$ that
\begin{align*}
\sum^{2^n} _{j=1} \Psi_2 \left(\frac{\ell(\Gamma_{n,j})^{2-\alpha}}{\lambda_2} \right)
\approx & \lambda^{\frac{1}{\alpha -1}} _2 \sum^{2^n} _{j=1} \ell(\Gamma_{n,j})^{\frac{2-\alpha}{1-\alpha}} \log^{\frac{\kappa \alpha}{\alpha -1}} \left(\frac{\lambda_2}{\ell(\Gamma_{n,j})^{2-\alpha}} \right) \notag \\
\lesssim & \lambda^{\frac{1}{\alpha -1}} _2 \sum^{2^n} _{j=1} \ell(\Gamma_{n,j})^{\frac{2-\alpha}{1-\alpha}} \log^{\frac{\kappa \alpha}{\alpha -1}} \left(\frac{1}{\ell(\Gamma_{n,j})} \right) \notag \\
\lesssim &\lambda^{\frac{1}{\alpha -1}} _2 \sum^{2^n} _{j=1} |\Delta_{n,j}| =1,
\end{align*}
where the first inequality comes from the fact that $\lambda_2 \ll 1$ when the generation $n$ of $\Gamma_{n,j}$ goes to infinity.
The second inequality is from \eqref{26.7.22_8}.
By the definition of the Luxemburg norm, we then obtain that 
\begin{equation}\label{26.7.22_9}
||\{\ell(\Gamma_{n,j})^{2-\alpha}\}||_{\Psi_2} \lesssim \lambda_2 = |\cup^{2^n} _{j=1} \Delta_{n,j}|^{1-\alpha} \le |\Omega|^{1-\alpha},
\end{equation}
where the equality holds because $\{\Delta_{n,j} : j=1,...,2^n \}$ are piecewise disjoint, and the rightmost inequality follows from the inclusion $\cup^{2^n} _{j=1} \Delta_{n,j} \subset \Omega .$
Therefore, we conclude from \eqref{26.7.22_10}, \eqref{26.7.22_7} and \eqref{26.7.22_9} that
\begin{align}
\sum^{2^n}_{j=1}\ell(\Gamma_{n,j})^{2-\alpha} \lesssim  & ||\{1\}||_{\Phi_2} \, ||\{\ell(\Gamma_{n,j})^{2-\alpha}\}||_{\Psi_2} \notag \\
\lesssim & 2^{\alpha n} \log^{\alpha \kappa } (2^n) |\Omega|^{1-\alpha} . \label{26.7.22_6}
\end{align}

Given $\alpha \in (0,1)$ as in \eqref{26.7.20_1}, we let $p=2-\alpha .$
It follows from \eqref{26.7.22_6} that 
\begin{align*}\label{26.7.22_5}
\sum^{\infty}_{n=n_0} 2^{(p-2)n} \sum^{2^n}_{j=1} \ell(\Gamma_{n,j})^p
= & \sum^{\infty}_{n=n_0} 2^{-\alpha n} \sum^{2^n}_{j=1}\ell(\Gamma_{n,j})^{2-\alpha} \\
\lesssim & \sum^{\infty}_{n=n_0} 2^{-\alpha n} 2^{\alpha n} \log^{\alpha \kappa } (2^n) |\Omega|^{1-\alpha}  \\
= & |\Omega|^{1-\alpha}  \sum^{\infty}_{n=n_0} \log^{\alpha \kappa } (2^n) <\infty,
\end{align*}
where the final inequality follows from $|\Omega|<\infty$ and $\alpha \kappa = - (\beta + \eta)<-1 .$
Given a homeomorphism $\varphi : \mbs \ra \pao ,$  we thus conclude from Theorem \ref{23.12.10_1} that there is a homeomorphic extension $h:\overline{\mbd}\to \overline{\Omega}$ that belongs to
$W^{1,2-\alpha}(\mbd,\Omega).$
Moreover, the Sobolev norm is bounded by a constant $M = M (\Omega) ,$ which is independent of $\varphi$ and $h .$ 
The proof is complete.
\end{proof}

\begin{proof}[Proof of Theorem \ref{26.7.20}]
The proof is analogous to the arguments for the construction of Example \ref{23.12.12_1}. 
We follow the notation introduced earlier and sketch the new arguments. 
The parameter settings need to be modified.

\bigskip
{\it Construction of $\Omega .$}
\bigskip

Denote by $R_n$ for $n=1,2,...$ open rectangles  with side lengths 
\begin{equation*}
a_n = 2^{-\frac{n}{2}} n^{-\epsilon} \qquad \mbox{ and } \qquad
b_n = 2^{-\frac{n}{2}} n^{-1-2\epsilon} ,
\end{equation*}
where $\epsilon >0$ is small enough.
As in Example \ref{23.12.12_1},
from $R_n $  we construct a snake-like open set $S_n .$ 
The length of the middle axis of $S_n$ is approximately 
\begin{equation}\label{23.7.24_1}
\ell_n = 2^{n \frac{1-\alpha}{\alpha-2}} n^{-1},
\end{equation}
where $\alpha \in (0,1) ,$
and the diameter of its opening is 
\begin{equation}\label{26.7.24_2}
o_n = \ell^{\frac{1}{1-\alpha}} _n \log^{\frac{\beta}{1-\alpha}} \left(\frac{1}{\ell_n} \right)
\end{equation}
where $\beta >2 .$
Denote by $d_n$ the distance between $S_n$ and  the $b_n$-sides of $R_n .$
On each $b_n$-side of $R_n ,$ we glue one $R_{n+1} .$
To connect $S_n$ and $S_{n+1} ,$ we use a rectangle $E_n$ whose side lengths are $o_n$ and $o_n + d_n + b_{n+1} .$
Let $T_n = S_n \cup E_n .$ 
Denote by $\Omega$ the union of all $T_n $ for all $ n \ge 1 .$
Let $\pao$ be the boundary of $\Omega .$
Analogously to Example \ref{23.12.12_1}, 
it follows from the inequality \eqref{25.2.18_1} that $R_{n-1}$ does not intersect $R_k$ for any $k \ge  n .$
Therefore, $\pao$ does not have self-intersections.

\bigskip
{\it Parametrization $\varphi$ of $\pao .$}
\bigskip

Analogously to the decomposition in Example \ref{23.12.12_1}, we 
decompose $\pao$ as
\begin{equation*}
\pao = \Gamma_1 \cup  \Gamma_2 \cup  \Gamma_3 \cup  \Gamma_4 .
\end{equation*}
Denote by $Q$ a square with the side length $\sum^{\infty}_{n=1} 2^{n-1} \hat{\ell}_n  ,$ where 
\begin{equation}\label{26.7.24_5}
\hat{\ell}_n = 2^{-n} n^{-2} .
\end{equation}  
Decompose the boundary of $Q$ as $\partial Q = I_1 \cup I_2 \cup  I_3 \cup I_4 .$
Analogously to Example \ref{23.12.12_1}, we parametrize $\pao$ by piecewise homeomorphisms 
\begin{equation}\label{26.7.24_4}
\varphi : I_i \ra \Gamma_i
\end{equation}
for $i=1,2,3,4 .$
Therefore, $\Omega$ is a Jordan domain.

\bigskip
{\it The estimate $\rho_{\Omega} (z_1 ,z) \lesssim d(z,\partial  \Omega)^{-\alpha} \log^{-\beta} (1/d(z,\partial  \Omega)) .$}
\bigskip

Notice that quasihyperbolic metric $h_{\Omega} (\cdot ,\cdot)$ is comparable to the hyperbolic metric $\rho_{\Omega} (\cdot ,\cdot)$ in a simply connected domain $\Omega .$ 
It then suffices to show the desired inequality for $h_{\Omega} (z_1 ,\cdot) .$
For each enlarged snake-like domain $T_n ,$ we denote by $z_n$ the middle point of the middle axis of $T_n .$
Fix $z \in \Omega.$ By the construction of $\Omega ,$ there is $T_n$ for some $n \ge 1 $ such that $z \in T_n  .$
Then $d (z , \partial T_n) = d(z ,\pao) .$
Let $\hat{z}$ be the projection of $z$ on the middle axis of $T_n .$
Denote by $\lambda (z_n, \hat{z})$ the  part of middle axis in $T_n$ joining $z_n$ to $\hat{z} .$
Then $\ell(\lambda (z_n, \hat{z})) \lesssim \ell_n ,$ where $\ell_n$ is from \eqref{23.7.24_1}.
We then calculate from the value of $o_n$ in \eqref{26.7.24_2} that 
\begin{align}
h_{\Omega} (z_n , \hat{z}) \le & \int_{\lambda (z_n, \hat{z})} \frac{d \xi}{d(\xi , \partial T_n)} 
=\frac{2 \ell(\lambda (z_n, \hat{z}))}{o_n} 
\lesssim   \frac{\ell_n}{o_n}  \notag\\
= & \ell^{\frac{\alpha}{\alpha -1}} _n \log^{\frac{\beta}{\alpha-1}} \left(\frac{1}{\ell_n} \right) .\label{26.7.23_1}
\end{align}
Again, by the value of $o_n$ we have
\begin{align}
d(\hat{z},\partial  \Omega)^{-\alpha} \log^{-\beta} \left(\frac{1}{d(\hat{z},\partial  \Omega)} \right)
 = & \left(\frac{o_n}{2} \right)^{-\alpha}  \log^{-\beta} \left(\frac{2}{o_n} \right)\notag \\
\approx & \ell^{\frac{\alpha}{\alpha -1}} _n  \log^{\frac{\beta}{\alpha -1}} \left(\frac{1}{\ell_n} \right). \label{26.7.23_2}
\end{align} 
Combining \eqref{26.7.23_1} with \eqref{26.7.23_2}, we arrive at
\begin{align}
h_{\Omega} (z_n , \hat{z}) \lesssim & d(\hat{z},\partial  \Omega)^{-\alpha} \log^{-\beta} \left(\frac{1}{d(\hat{z},\partial  \Omega)} \right) \notag \\
\le & d(z,\partial  \Omega)^{-\alpha} \log^{-\beta} \left(\frac{1}{d(z,\partial  \Omega)} \right) , \label{26.7.23_3}
\end{align}
where the second inequality follows from the monotonicity of $f(t) = t^{-\alpha} \log^{-\beta} (1/t)$ and $d(\hat{z},\partial  \Omega) \ge  d(z,\partial  \Omega) .$
Denote $\gamma(\hat{z} ,z)  $ the line segment joining $\hat{z}$ to $z .$
We estimate
\begin{align}
h_{\Omega} (\hat{z}, z) \le  & \int_{\gamma(\hat{z} ,z)  } \frac{d \xi}{d(\xi ,\partial T_n)} = \ln \frac{o_n }{2 d (z ,\partial T_n )} \notag  \\
\lesssim  & d(z,\partial  \Omega)^{-\alpha} \log^{-\beta} \left(\frac{1}{d(z,\partial  \Omega)} \right) . \label{26.7.23_4}
\end{align}
It follows from \eqref{26.7.23_3} and \eqref{26.7.23_4} that 
\begin{equation}\label{26.8.3_1}
h_{\Omega } (z_n, z) \le   h_{\Omega  } (z_n, \hat{z}) + h_{\Omega  } (\hat{z}, z) \lesssim d(z,\partial  \Omega)^{-\alpha} \log^{-\beta} \left(\frac{1}{d(z,\partial  \Omega)} \right) .
\end{equation}
Analogously to \eqref{26.7.23_1}, we have $h_{\Omega} (z_k ,z_{k+1}) \lesssim \ell^{\alpha / (\alpha -1)} _k \log^{\beta/(\alpha-1)} (1/\ell_k ) $ for any $k \ge 1 .$
Therefore, it follows from values of $\ell_n$ and $o_n$ that
\begin{align}
\sum^{n-1} _{k=1} h_{\Omega} (z_k, z_{k+1})  \lesssim & \sum^{n-1} _{k=1} \ell^{\frac{\alpha}{\alpha -1}} _k \log^{\frac{\beta}{\alpha-1}} \left(\frac{1}{\ell_k} \right)
 \notag \\
 \approx & \ell^{\frac{\alpha}{\alpha -1}} _n \log^{\frac{\beta}{\alpha-1}} \left(\frac{1}{\ell_n} \right)  
\approx  o^{-\alpha} _n \log^{-\beta} \left(\frac{1}{o_n} \right) . \label{26.7.23_5}
\end{align}
By the triangle inequality of $h_{\Omega} ,$ we conclude from  \eqref{26.7.23_5} and \eqref{26.8.3_1} that 
\begin{align*}
 h_{\Omega} (z_1, z) 
\le & \sum^{n-1} _{k=1} h_{\Omega} (z_k, z_{k+1}) +h_{\Omega } (z_n ,z) \\
\lesssim &  o^{-\alpha} _n \log^{-\beta} \left(\frac{1}{o_n} \right) +  d(z,\partial  \Omega)^{-\alpha} \log^{-\beta} \left(\frac{1}{d(z,\partial  \Omega)} \right) \\
\approx & d(\hat{z},\partial  \Omega)^{-\alpha} \log^{-\beta} \left(\frac{1}{d(\hat{z},\partial  \Omega)} \right) + d(z,\partial  \Omega)^{-\alpha} \log^{-\beta} \left(\frac{1}{d(z,\partial  \Omega)} \right) \\
\lesssim & d(z,\partial  \Omega)^{-\alpha} \log^{-\beta} \left(\frac{1}{d(z,\partial  \Omega)} \right) ,
\end{align*}
where the final inequality follows from the arguments for  \eqref{26.7.23_3}.
The desired inequality for  $h_{\Omega} (z_1 ,\cdot) $ is obtained.

\bigskip
{\it The positive part of Theorem \ref{26.7.20}}
\bigskip

Given $1 <p <2 ,$ we take $\alpha \in (0,1)$ such that $p=2-\alpha .$
Take the other parameter $\beta >2 .$
By  the above arguments,  there is a Jordan domain $\Omega$
such that the hyperbolic distance $\rho_{\Omega}$ satisfies the growth condition
\begin{equation*}
\rho_{\Omega}(z,z_0)\le C d(z,\partial  \Omega)^{-\alpha} \log^{-\beta} \left(\frac{1}{d(z,\partial  \Omega)} \right)
\end{equation*}
for some fixed $z_0\in \Omega$ and all $z\in \Omega ,$ where $C$  is a constant.
Denote by $\gamma_p$ the boundary of  $\Omega .$
Therefore, by Lemma \ref{26.7.24_3} there is a constant $M$ such that 
every parametrization of $\gamma_p$ via $\mbs \subset \mbs^2$ has a homeomorphic $W^{1,p}-$extension with norm bounded by $M.$

\bigskip
{\it  The negative part of Theorem \ref{26.7.20}}
\bigskip

Let $\varphi $ be the parametrization of $\gamma_p$ as in \eqref{26.7.24_4}.
For simplicity, we define this $\varphi$ on the boundary of a square $Q .$
Suppose that  $h : Q \ra \Omega$ is a homeomorphic Sobolev extension of $\varphi .$
Recall that $I_3$ is  a part of $\partial Q $ in Example \ref{23.12.12_1}.
We proceed to establish a lower bound for  the integral of $|Dh|^q,$ $q>1,$ over a sub-domain of $Q$ close to  $I_3 .$ 
As in \eqref{25.2.26_1}, we express $I_3$ as
\begin{equation}\label{26.7.24_7}
I_3 = \mathcal{C} \cup \cup^{\infty} _{i=1} \cup_{\lambda \in K_i} \lambda ,
\end{equation}
where $\mathcal{C} $ is a Cantor-type set and $\lambda$ is an open interval.
Fix an open interval $\lambda _n \in K_n$ and denote by $p_n$ the middle point of $\lambda_n .$ Let  $\lambda_{n-1} \in K_{n-1}$ be the one in $K_{n-1}$ nearest to $\lambda_n ,$ and  choose $\lambda_{n+1}$ in $K_{n+1} $  analogously.
Let $r_n$ be the distance between $p_n$ and the middle point of $\lambda_{n+1} ,$ and let $\hat{r}_n$ be the distance between $p_n$ and $\lambda_{n-1} .$ Set $A_n$ to be the intersection of $Q$ and the open annulus centered at $p_n $ with radii $r_n$ and $\hat{r}_n .$ 
By the value of $\hat{\ell}_n $ as in \eqref{26.7.24_5},  we estimate
\begin{equation}\label{26.7.24_8}
r_n = \frac{\hat{\ell}_n}{2} + \frac{\hat{\ell}_{n+1}}{2}+\sum^{\infty} _{k=n+2} \hat{\ell}_k 2^{k-n-2}  \approx 2^{-n} n^{-1}
\end{equation}
and 
\begin{equation*}
\hat{r}_n =   \frac{\hat{\ell}_n}{2} +\hat{\ell}_{n+1}  + 2\sum^{\infty} _{k=n+2} \hat{\ell}_k 2^{k-n-2} \approx 2 r_n .
\end{equation*}
Let $\sigma_n $ be the maximal square embedded in $A_n$ such that $A_n \setminus \sigma_n$ is symmetric.
Recall that $1 <p <2 $ in the statement of Theorem \ref{26.7.20}. 
Analogously to \eqref{24.5.16_5},  we have 
\begin{equation}\label{26.7.24_5_1}
\ell^p _n r^{-p} _n \lesssim M_{|Dh|^p } (z) 
\end{equation}
for all $z \in \sigma_n ,$ where $M_{|Dh|^p }$ is the maximal function of the zero extension of  $|Dh|^p $ from $Q$ to the entire plane.
Take $q >p .$
We raise   \eqref{26.7.24_5_1} to the power $q/p$ and integrate over $z \in \sigma_n$ to obtain 
\begin{equation}\label{26.7.24_6}
\ell^{q} _n r^{2-q} _n \approx \int_{\sigma_n} (\ell^p _n r^{-p} _n )^{\frac{q}{p}} dz \lesssim \int_{\sigma_n} M_{|Dh|^p } (z)^{\frac{q}{p}} dz .
\end{equation}
Let $\Sigma_n$ be the collection of all squares $\sigma_n$ corresponding to $\lambda_n \in K_n ,$ where $K_n$ are as in \eqref{26.7.24_7}. 
For any $\sigma_i \in \Sigma_i$ and $\hat{\sigma}_j \in \Sigma_j ,$
we claim that $\sigma_i \cap  \hat{\sigma}_j = \emptyset .$
It suffices to check that $ A_i \cap  \hat{A}_j = \emptyset ,$ where $ \hat{A}_j $ is the half-annulus corresponding to $\hat{\sigma}_j .$
In fact, if $i=j $ we can show $A_i \cap  \hat{A}_j = \emptyset $ by the definition of $\hat{r}_i .$ If $i \neq j ,$ it suffices to consider the case $j =i+1 .$ 
It follows from $\hat{r}_{i+1} <r_i $ that $ A_i \cap \hat{A}_{i+1} = \emptyset .$
Since $\Sigma_n$ consists of  $2^{n-1}$ elements, we obtain from \eqref{26.7.24_6} that
\begin{align}
\sum^{\infty} _{n=1} 2^{n-1} \ell^{q} _n r^{2-q} _n 
 \lesssim & \sum^{\infty} _{n=1} \sum_{\sigma \in \Sigma_n}\int_{\sigma} M_{|Dh|^p } (z)^{\frac{q}{p}} dz  \notag\\
 \le & \int_{\mbr^2} M_{|Dh|^p } (z)^{\frac{q}{p}} dz 
\lesssim  \int_{Q} |Dh|^{q}  (z) dz, \label{26.7.24_9}
\end{align}
where the final inequality is from the $L^{q/p} (\mbr^2)$-boundedness of the maximal operator.
By inserting the values of  $\ell_n$  and $r_n $ from \eqref{23.7.24_1}  and \eqref{26.7.24_8},  we calculate that
\begin{equation*}\label{26.8.3_2}
\sum^{\infty} _{n=1} 2^{n-1} \ell^{q} _n r^{2-q} _n 
\approx \sum^{\infty} _{n=1} 2^{n (\frac{q}{2-\alpha}-1)}   n^{-2} = \infty
\end{equation*}
whenever $q >p ,$  where $p =2-\alpha .$
It follows from \eqref{26.7.24_9} that $h$ is not in $W^{1,q} (Q,\Omega) $ for any $q> p .$
\end{proof}

\section{Proof of Theorem \ref{posa}}\label{26.4.10_1}

\begin{proof}[Proof of Theorem \ref{posa}]
Let $Q = [0,1] \times [0,1]$ be the closed unit square. Denote by $u$ and $v$ the upper-left and lower-right vertices of $Q .$ 
We will construct a suitable Jordan arc $\Gamma \subset Q$ of  positive two-dimensional Lebesgue measure so that it joins $u$ to $v. $
The desired Jordan curve is then the union of $\Gamma$ and  $\mbs \setminus Q .$ 
By the Jordan curve theorem, there are two components in $\mbr^2 \setminus \Gamma :$ the bounded one $\Omega$ and the unbounded one $\Omega^* .$
Then $\pao=\Gamma$ is their common boundary.
We will show that  there are $C_1 >0$ and $z_0 \in \Omega$ such that 
\begin{equation}\label{25.2.5_3}
h_{\Omega} (z_0 ,z) \le C_1  \log ^3 \left( \frac{1}{d (z ,\pao)} \right)
\end{equation}
for all $z \in \Omega $ and there are $C_2 >0$ and $w_0 \in \Omega^*$ such that 
\begin{equation}\label{25.2.7_1}
h_{\Omega^* } (w_0 , w) \le C_2 \log ^3 \left( \frac{1}{d (w ,\pao)} \right)
\end{equation}  
for all $w \in \Omega^* .$

\bigskip
{\it Construction of the Jordan arc $\Gamma .$}
\bigskip

For $n= 1,2,...,$ we take 
\begin{equation}\label{25.2.5_1}
d_n = 3^{-n} \ln^{-2} (3^n) \qquad \mbox{ and } \qquad  \ell_{n+1} = 3^{-n} (1- 2 \sum^n _{k=1} 3^{k-1} d_k) .
\end{equation}
Set $\ell_1 = 1 .$ Our construction is based on pairs $R_n$ of slim horizontal rectangles and pairs $\hat{R}_n$ of slim vertical rectangles of side lengths $\ell_{n}$ and $d_n .$  
The distance between the rectangles in the same pair will be $\ell_{n+1}.$ Write $\partial R_n$ (or $\partial \hat{R}_n$) for the union of the boundaries of the two rectangles in $R_n$ (or in $\hat{R}_n$, respectively).


{\bf The first step.} We remove $R_1$ and $\hat{R}_1$ from $Q$ so that  $Q \setminus (R_1 \cup \hat{R}_1)$ consists of $9$ equivalent closed squares.
By a snake-shaped order, we number them as $Q^1 _1 , Q^2 _1 ,..., Q^9 _1  $ so that $u$ is in the `head" $ Q^1 _1$ and $v$ in the ``tail" $ Q^9 _1 .$
Join $Q^j _1 $ to $Q^{j+1} _1$ by an open line segment $\gamma^j _1 \subset \partial R_1 \cup \partial \hat{R}_1$ of length $d_1 .$ 
This results in a snake-shaped set $T_1,$ built from $Q^j _1 , \gamma^j _1 ,  Q^{j+1} _1$ in order for $j =1,...,8 .$
This set is shown on the left  of Figure \ref{25.2.15_1}.

\begin{figure}[htbp]
\centering
\begin{subfigure}[b]{.49\textwidth}
  \centering
  \includegraphics[height = 0.25\textheight]{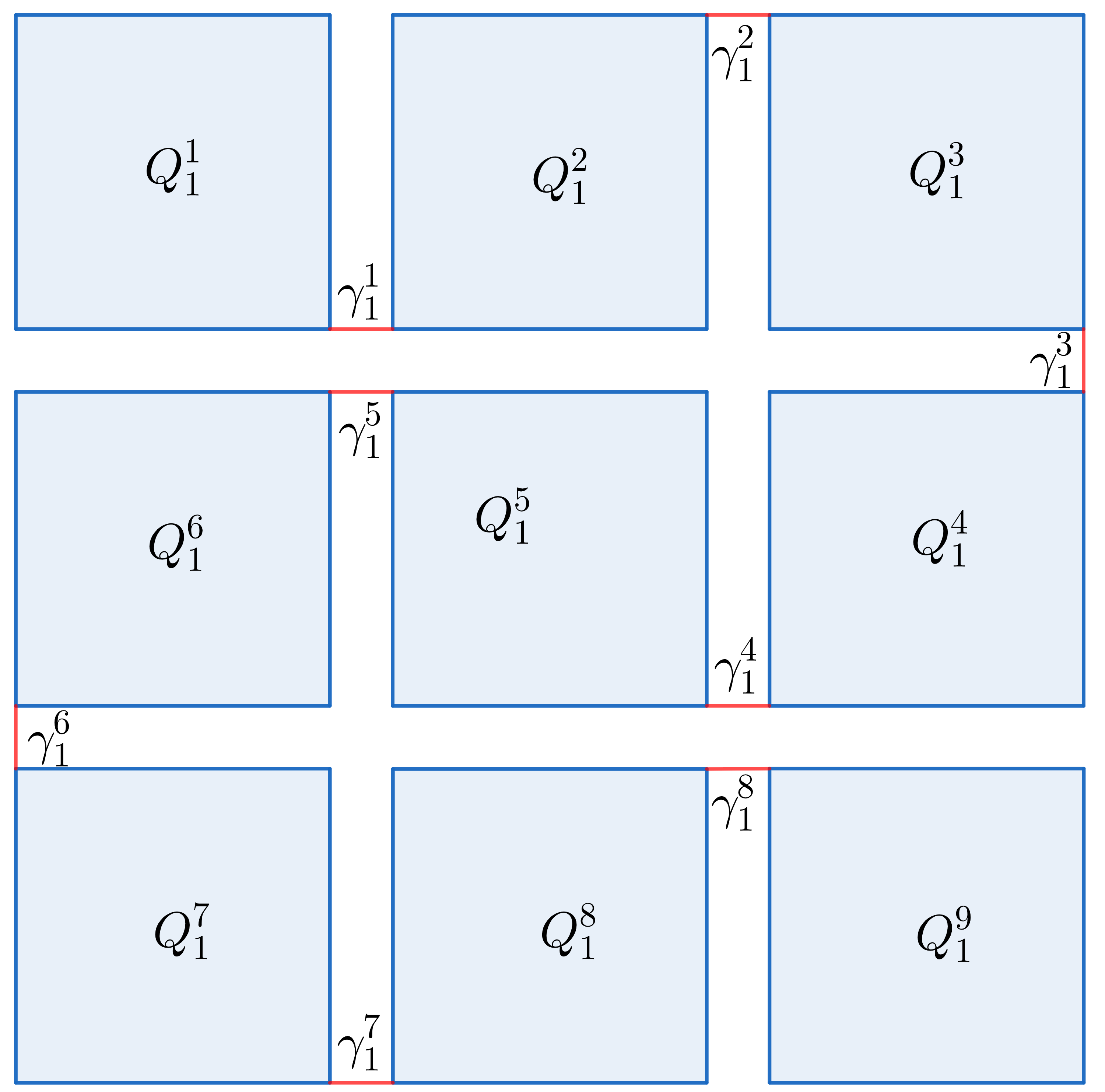}  
  \caption*{The $T _1$ set}
\end{subfigure}
\begin{subfigure}[b]{.49\textwidth}
  \centering
  \includegraphics[height = 0.248\textheight]{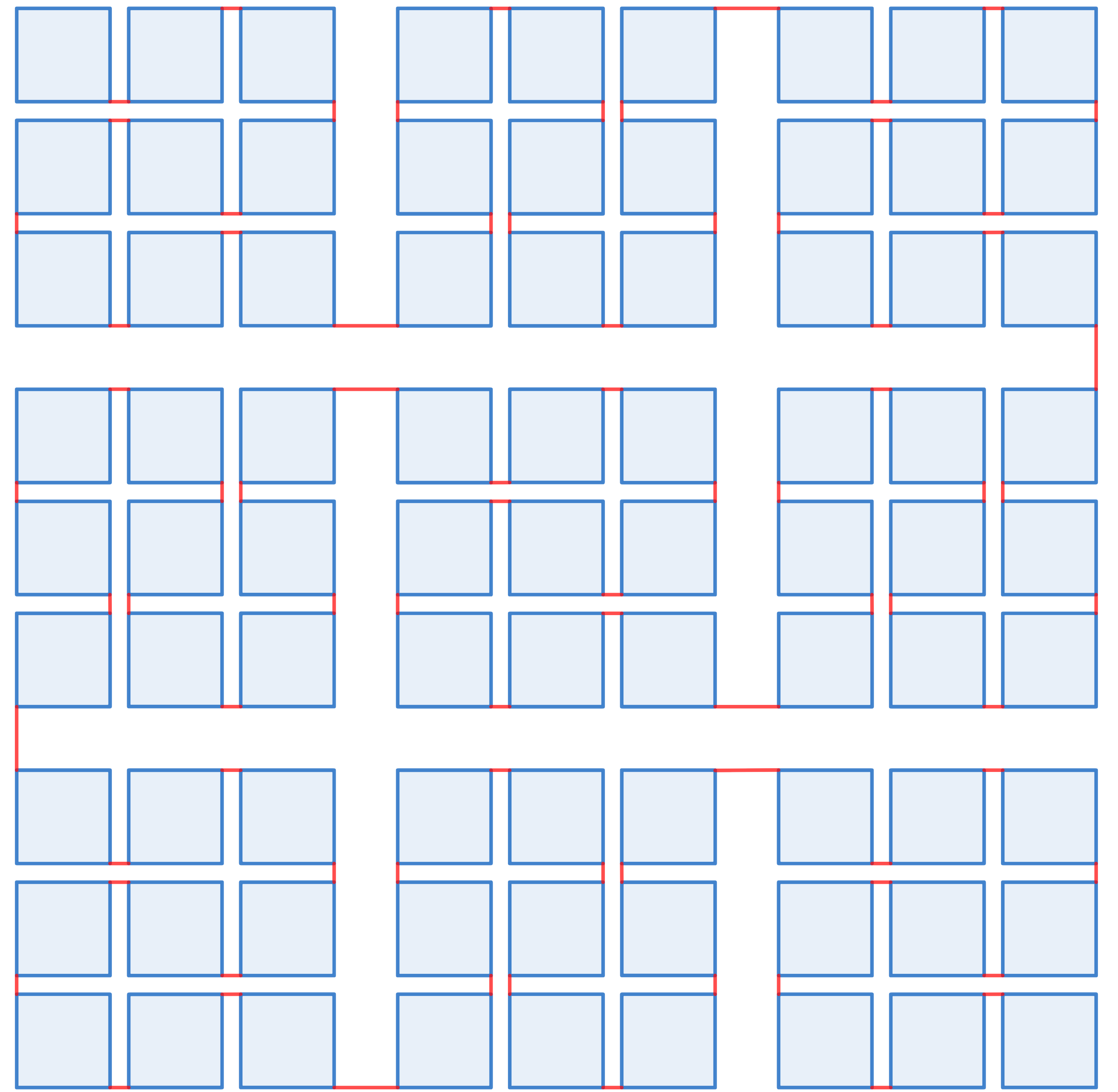}  
  \caption*{The $T _2$ set}
\end{subfigure}
\caption{The $T_1$ and $T_2$ sets}
\label{25.2.15_1}
\end{figure}

{\bf The second step.} 
We will construct a new snake-shaped set $T _2 $ joining $u$ to $v .$
This set is shown on the right of Figure \ref{25.2.15_1}.
To do this, we replace each square in $T_1$ with a variant of $T_1 .$
There are four variants
\begin{equation*}
T^1 _2 , \  T^2 _2 , \ T^3 _2 \mbox{ and } T^4 _2 .
\end{equation*}
See Figure \ref{25.2.14_2}.
We first construct $T^1 _2 $  and replace $Q^1 _1 $ with it.
Denote by $u^1 _1$ and $v^1 _1$ the upper-left vertex and lower-right vertex of $Q^1 _1 .$
From $Q^1 _1$ we delete $R_2$ and $\hat{R}_2$ so that $Q^1 _1 \setminus (R_2 \cup \hat{R}_2 )$ consists of $9$ equivalent closed squares.
Number them as $Q^j _2$ for $j=1,...,9 ,$ and join $ Q^j _2$ to $Q^{j+1} _2$ by an open line segment $\gamma^j _2 \subset \partial R_2 \cup \partial \hat{R}_2$ of length $d_2 .$ 
Then we obtain a snake-shaped set $T^1 _2 $ joining $u^1 _1$ to $v^1 _1 .$
This $T^1 _2 $ is a scaled-down copy of  $T_1 $ with $Q^j _2 , \gamma^j _2 ,  Q^{j+1} _2$ in order for $j =1,...,8 .$
Moreover, $u^1 _1$ is the entry point of $T^1 _2 $ and $v^1 _1$ is its exit point.
We rotate $T^1 _2 $ counterclockwise around  the center of $Q^1 _1 $ by degrees $\pi /2 , \pi$ and $3 \pi /2$ to obtain $ T^2 _2 , T^3 _2 $ and $T^4 _2. $ 
\begin{figure}[htbp]
\centering
\begin{subfigure}[b]{.2\textwidth}
  \centering
  \includegraphics[width=\linewidth]{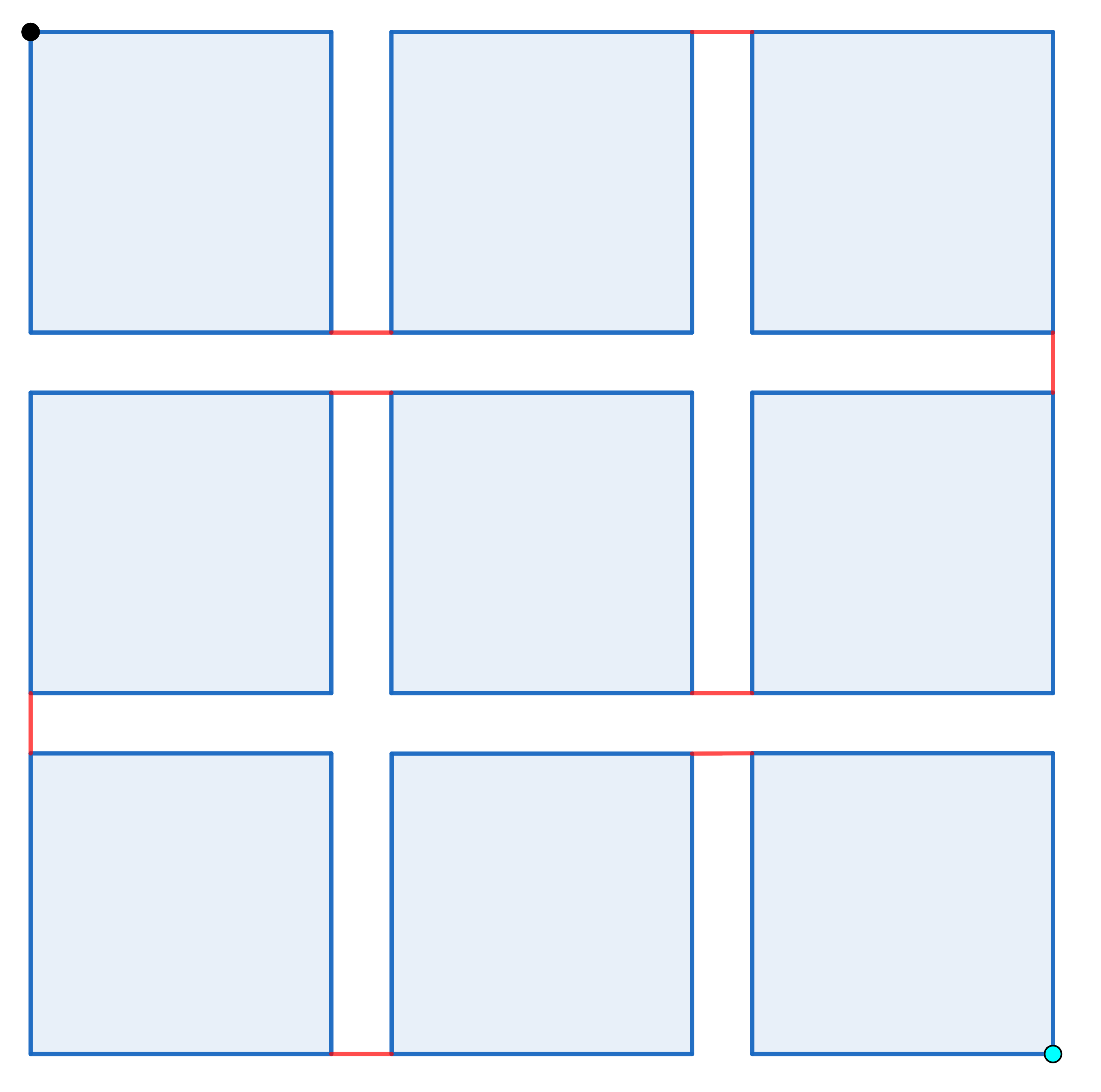}  
  \caption*{$T^1 _n$}
\end{subfigure}
\begin{subfigure}[b]{.2\textwidth}
  \centering
  \includegraphics[width=\linewidth, angle=90 ]{25.2.14_2}   
  \caption*{$T^2 _n$}
\end{subfigure}
\begin{subfigure}[b]{.2\textwidth}
  \centering
  \includegraphics[width=\linewidth, angle=180]{25.2.14_2}   
  \caption*{$T^3 _n$}
\end{subfigure}
\begin{subfigure}[b]{.2\textwidth}
  \centering
  \includegraphics[width=\linewidth, angle=270]{25.2.14_2}   
  \caption*{$T^4 _n$}
\end{subfigure}

\caption{Snake-shaped sets $T^1 _n ,\ T^2 _n ,\ T^3 _n ,$ and $T^4 _n $ for $n \ge 2 .$}
\label{25.2.14_2}
\end{figure}
Then, the beginning points of $ T^2 _2 , T^3 _2 ,$ and $T^4 _2 $ are at the lower-left, the lower-right, and the upper-right vertices of $Q^1 _1 .$ The corresponding exit point is the other endpoint of the diagonal line. We refer to different $T^i _2$ by the location of their beginning points.

We next replace $Q^j _1$ for $j =2,...,9$ with $T^i _2 .$
Assume that one $Q^{j} _1$ has been replaced by a $T^i _2 .$
Denote by $v^{j} _1$ the ending point of this  $T^i _2 .$
Then we know the location of $v^{j} _1$  in $Q^{j} _1 . $
Moreover, $v^{j} _1$ is the joining point between $Q^{j} _1$ and  $\gamma^{j} _1 .$
Denote by $u^{j +1} _1$ the other endpoint of $\gamma^{j} _1 .$
The location of $u^{j +1} _1$ in $Q^{j +1} _1$ determines the beginning point of a new $T^{i'} _2 .$
Replace $Q^{j+1} _1 $ with this $T^{i'} _2 .$
Continue this process until we replace the last $Q^9 _1$ with a $T^i _2 .$
This gives rise to a snake-shaped set $T_2$ joining  $u$ to $v .$
It is composed of $Q^j _2 ,\ \gamma_2 ^{k_2}$ and $\gamma_1 ^{k_1}$ for $j=1,...,9^2 ,$ $k_2 = 1,...,72$ and $k_1 = 1,...,8 .$

{\bf A general step.} After the step $n$, we have obtained a snake-shaped set $T_n$ joining $u$ to $v .$
This $T_n$ is composed of closed squares $Q^j _n$ and open line segments $\gamma^{k_{\tau}} _{\tau}$ for $j =1,...,9^n  ,$ $\tau = 1,..., n$ and $k_{\tau} = 1,...,9^{\tau -1} 8 .$ 
Moreover, $u$ is the entry point of $T_n$ and $v$ is its exit point.
Along this direction,  $Q^j _n $ and $\gamma^{k_{\tau}} _{\tau} $ alternate in $T_n .$
Fix $\tau = 2,...,n .$ 
The eight $\gamma^{k_{\tau}} _{\tau}$ form a group in this direction. After such a group, there comes $\gamma^{k_{\tau'}} _{\tau'}$ with $\tau' < \tau .$

We will construct a new snake-shaped set $T_{n+1}$ joining $u$ to $v .$ 
To do this, we replace each $Q^j _n$ in $T_n$ with a variant of $T_1 .$
There are four variants 
$$T^1 _{n+1} , T^2 _{n+1} , T^3 _{n+1}  \mbox{ and } T^4 _{n+1} .$$ See Figure \ref{25.2.14_2}.
We first construct $T^1 _{n+1}$ and replace $Q^{1} _n$ by it.
Denote by $u^1 _n$ and $v^1 _n$ the upper-left  and lower-right vertices of $Q^1 _1 .$
We remove $R_{n+1}$ and $\hat{R}_{n+1}$ from $Q^1 _n$ so that $Q^1 _n \setminus (R_{n+1} \cup \hat{R}_{n+1} )$ consists of nine equivalent closed squares. 
Number them as $Q^j _{n+1}$ for $j =1,...,9$ and connect $ Q^j _{n+1}$ and $Q^{j+1} _{n+1}$ by an open line segment $\gamma^j _{n+1} \subset \partial R_{n+1} \cup \partial \hat{R}_{n+1}$ of length $d_{n+1} .$ 
Then we obtain a snake-shaped set $T^1 _{n+1} $ composed of $Q^j _{n+1} , \gamma^j _{n+1}  ,  Q^{j+1} _{n+1} $ in order for $j =1,...,8 .$
This $T^1 _{n+1} $ is a scaled-down copy of $T_1 $ with the entry point $u^1 _n$ and the exit point $v^1 _n .$
Rotate $T^{1} _{n+1}$ counterclockwise by degrees $\pi /2 , \pi$ and $3 \pi /2$ to obtain $ T^2 _{n+1} , T^3 _{n+1} $ and $T^4 _{n+1} ,$ respectively. 
Then, the beginning points of $ T^2 _{n+1} , T^3 _{n+1}  $ and $T^4 _{n+1} $ are at the lower-left, the lower-right and the upper-right vertices of $Q^1 _n .$ The corresponding ending point is the other endpoint of the diagonal line in $Q^1 _n .$
We refer to different $T^i _{n+1}$ by the location of their beginning points.

We next replace $Q^j _n$ for $j =2,...,9^n$ with $T^i _{n+1} .$
This process is inductive.
Assume that $Q^j _n$ has been replaced by  $T^i _{n+1} .$
Denote by $v^j _n$ the endpoint of this $T^i _{n+1} .$
 According to the construction of $T_n ,$ this  $v^j _n$ is the joining point between $Q^j _n$ and a line segment $\gamma^{k_{\tau}} _{\tau} .$
Denote by $u^{j+1} _n$ the other endpoint of $\gamma^{k_{\tau}} _{\tau} .$ 
Then $u^{j+1} _n$  is a vertex of the next square $Q^{j+1} _n .$
The location of $u^{j+1} _n$ in $Q^{j+1} _n $ determines the beginning point of the new $T^{i'} _{n+1} .$
We replace $Q^{j+1} _n $ with this $T^{i'} _{n+1} $
and continue this process until we have replacec all $Q^j _n $ for $j =2,...,9^n .$

By induction, we obtain countably many open segments $\gamma^{k_{n}} _{n}$ for $n \in \mbn$ and $k_{n} = 1,2, ...,9^{n-1} 8  .$ From the complement of their union, $( \cup^{\infty} _{n=1} \cup^{9^{n-1} 8 } _{k_{n} =1 }  \gamma^{k_{n}} _{n})^c$ ,  we define a set $\mcc $  as the collection of accumulation points of $ \cup^{\infty} _{n=1} \cup^{9^{n-1} 8 } _{k_{n} =1 }  \gamma^{k_{n}} _{n} .$ 
Define $\Gamma$ as
\begin{equation}\label{24.6.20_1}
\Gamma = \mcc \cup  \cup^{\infty} _{n=1} \cup^{9^{n-1} 8 } _{k_{n} =1 }  \gamma^{k_{n}} _{n} .
\end{equation}
To show $\Gamma$ is a Jordan arc, we will construct a parametrization of $\Gamma .$

\bigskip
{\it Parametrization of  $\Gamma .$ }
\bigskip

Let $I=[0,1] $ be the closed unit interval. 
We decompose $I$ via a Cantor set construction. 
For $n \in \mbn ,$ we set
\begin{equation*}
\delta_n = 17^{-n} .
\end{equation*}
Step $1:$  We remove an open set $K_1$ composed of eight open intervals of equal length from $I$ so that $I_1 = I \setminus K_1 $ 
consists of nine closed disjoint intervals, all of length $\delta_1.$
In the $2$nd step, we remove from every component of $I_1$ an open set $K_2$ consisting of  eight distinct open intervals of equal length so that
 $I_2 = I_1 \setminus K_2$  is composed of  $9^2$  distinct closed intervals of length $\delta_2.$
In general, after step $n$, we have $9^n$ closed intervals of length $\delta_n.$ Denote by $I_n$ the union of them. 
In step  $n+1$, we delete from each component of $I_n$ eight distinct open intervals of equal length so that the remaining closed intervals all have length $\delta_{n+1}.$
Let $K_{n+1}$ be the collection of all open intervals in this step and set $I_{n+1} = I_n \setminus K_{n+1} .$ 
Then $K_{n+1}$ consists of $9^{n} 8 $ open intervals and $I_{n+1}$ of $9^{n+1}$ of closed interval. 
Continuing by induction, we obtain a decomposition 
\begin{equation}\label{24.6.20_2}
I = E \cup \cup^{\infty} _{n=1} \cup_{\lambda \in K_n} \lambda ,
\end{equation}
where $E = \cap^{\infty} _{n=1} I_n$ is a Cantor-type set.

We next construct a parametrization of $\varphi : I \ra \Gamma$ via \eqref{24.6.20_1} and \eqref{24.6.20_2}.
For each $n \in \mbn ,$ the set $K_n$ consists of $9^{n-1} 8$  open line segments and so does the set $\cup^{9^{n-1} 8 } _{k_{n} =1 }  \gamma^{k_{n}} _{n} .$ 
The snake-shaped set $T_n $ gives a natural order for  $\gamma^{1} _{n} ,...,  \gamma^{9^{n-1} 8 } _{n} .$
 We consider the increasing order for $K_n$ in  $.I $ 
We map the  components of  $K_n$ to the corresponding $\gamma^{k_{n}} _{n}$ by a piecewise affine mapping $\varphi .$ 
Denote by $ \partial K_n$ the collection of endpoints of the intervals in $ K_n $ and by $\overline{K_n}$ the collection of the closure of intervals in $ K_n .$
Thus, $\varphi$ extends to $\overline{K_n}$  as a homeomorphism.
It remains to define $\varphi$ on $E \setminus \cup^{\infty} _{n=1} \partial K_n .$ 
For $z $ in this set, we take $z_n \in \partial K_n$ such that $|z-z_n| = \inf_{\xi \in \cup_{\lambda \in K_n} \lambda} |z-\xi | .$ Then $\lim_{n \ra \infty} z_n =z .$
By the construction of $\Gamma ,$ we have
\begin{equation*}
|\varphi (z_n) - \varphi (z_{n+1})| \le \ell_{n+1} \le \frac{1}{3^n} .
\end{equation*}
Then $\{\varphi (z_n)\}$ is a Cauchy sequence and $\lim_{n \ra \infty} \varphi (z_n ) $ exists. In summary, we define
\begin{equation}\label{25.3.6_1}
\varphi (z) = 
\begin{cases}
\mbox{piecewise affine}, & \mbox{ if } z \in \overline{K_n} \\
\lim_{n \ra \infty} \varphi (z_n ),  & \mbox{ if } z  \in I \setminus \cup^{\infty} _{n=1} \overline{K_n } .
\end{cases}
\end{equation}
Obviously $\varphi (I) \subset \Gamma .$ To show that $\varphi : I \ra \Gamma$ is a homeomorphism, it suffices to prove that $\varphi$ is continuous and bijective.

To establish the continuity of $\varphi ,$ it suffices to construct a sequence of continuous $\varphi_m$ on $I$ converging uniformly to $\varphi .$
The snake-shaped set $T_m$ contains the squares $Q^j _m$ and the line segments $\gamma^{k_{\tau}} _{\tau}$ for $j =1,...,9^m  ,$ $\tau = 1,..., m$ and $k_{\tau} = 1,...,9^{\tau -1} 8 .$ 
Replace each $Q^j _m$  with a diagonal $L^j _m$ of $Q^j _m $
so that  $L^j _m$ connects the entry point into $Q^j _m$ to the exit point. By the decomposition \eqref{24.6.20_2}, we have
$I = I_m \cup \cup^m _{n=1} \cup_{\lambda \in K_n} \lambda .$
We map the interval in $I_m $  onto the  corresponding  $L^j _m$ by a piecewise affine mapping $\psi_m.$
Let  
\begin{equation*}
\varphi_m (z) = 
\begin{cases}
\varphi(z), & \mbox{ if } z \in \overline{K_n}  \mbox{ with } 1 \le n \le m ,\\ 
\psi_m, & \mbox{ if } z \in I_m .
\end{cases}
\end{equation*}
Then $\varphi_m$ is a homeomorphism on $I .$
Note that $\varphi_m (z) = \varphi (z)$ for $z \in \cup^{m} _{n=1} \overline{K_n} .$
Given $z \in I \setminus \cup^{m} _{n=1} \overline{K_n} ,$ we 
\begin{equation*}
|\varphi(z) - \varphi_m (z)| \le \diam\, Q^j _m  .
\end{equation*}
Note that $\diam\, Q^j _m  \approx \ell_m \approx 3^{-m}$ by \eqref{25.2.5_1}.
Hence, $\varphi_m$ converges uniformly to $\varphi$ over $I .$

To see that $\varphi$ is surjective, let $w \in \Gamma \setminus \cup^{\infty} _{n=1} \cup^{9^{n-1} 8 } _{k_{n} =1 }  \overline{\gamma^{k_{n}} _{n}} .$
There exists $w_n \in \cup^{9^{n-1} 8 } _{k_{n} =1 }  \partial \gamma^{k_{n}} _{n} $ such that $\lim_{n \ra \infty} w_n =w .$ Take $z_n \in \partial K_n$ such that $w_n =\varphi (z_n) .$ Since $\{z_n\}$ is bounded, there is a subsequence, denoted by $\{\hat{z}_n\}$, converging to a point $z \in I .$ Hence, by the continuity of $\varphi$ we obtain
\begin{equation*}
\varphi(z) = \lim_{n \ra \infty} \varphi (\hat{z}_n) = \lim_{n \ra \infty} w_n =w .
\end{equation*}
Therefore, $\varphi$ is surjective.

Fix distinct points $z ,\ \hat{z} \in I .$ We show that  $\varphi (z ) \neq \varphi (\hat{z} ) .$ If $z , \ \hat{z} \in I \setminus E ,$ then $\varphi (z ) \neq \varphi (\hat{z} ) $ is obvious by \eqref{25.3.6_1} and the geometry of $\Gamma .$
Assume that $z \in E$ and $\hat{z} \in I \setminus E. $
By \eqref{24.6.20_2}, there is $n_0 \ge 1$ such that $\hat{z}$ is in a component of $K_{n_0} .$ Moreover, by the definition of $\varphi$ there is an open interval $\gamma^{k_n} _n $ with $\varphi(\hat{z}) \in \gamma^{k_n} _n  .$  
Since $\varphi (E) \subset \mathcal{C}$ and $\mathcal{C} \cap (\cup^{\infty} _{n=1} \cup^{9^{n-1} 8 } _{k_{n} =1 }  \gamma^{k_{n}} _{n})  = \emptyset ,$ we have $\varphi (z ) \neq \varphi (\hat{z} ) .$ 
Assume that $z , \ \hat{z} \in   E $ with $z < \hat{z} .$
Take $n_0$ as the minimal integer such that
there is a component $\tau_{n_0}$ in $K_{n_0}$ with $\tau_{n_0} \subset (z, \hat{z}) .$
By the decomposition \eqref{24.6.20_2}, we have $I= I_{n_0} \cup \cup^{n_0} _{i=1} \cup_{\lambda \in K_i} \lambda $ and $z , \ \hat{z} \in  \cup_{\tau \in I_{n_0}} \tau .$ It follows from \eqref{25.3.6_1} that $\varphi (I_{n_0}) \subset \cup^{9^{n_0}} _{j=1} Q^j _{n_0} $ and 
$\varphi (I) \subset T_{n_0} .$
Hence, there are $Q^j _{n_0} $ and $Q^k _{n_0}$ with $j <k$ such that $\varphi (z ) \in Q^j _{n_0}$ and $\varphi (\hat{z} ) \in Q^k _{n_0}.$
Moreover, $Q^j _{n_0} ,\ \tau_{n_0} $ and $\ Q^k _{n_0} $
appear in $T_{n_0} $ in order.
By the geometry of $T_{n_0} ,$ the distance between $Q^j _{n_0}$ and $Q^k _{n_0}$ is larger than $d_{n_0} .$ Then $|\varphi (z )  - \varphi (\hat{z} )| \ge d_{n_0} .$ In conclusion, $\varphi$ is injective on $I .$

\bigskip
{\it Positivity  of the area of the arc $\Gamma$  }
\bigskip

When constructing the arc $\Gamma $ in \eqref{24.6.20_1}, we introduced a sequence of snake-shaped sets $\{T_n \}.$ Every $T_n$ is composed of squares $Q^j _n$ and line segments $\gamma^{k_{\tau}} _{\tau}$ for $j =1,...,9^n  ,$ $\tau = 1,..., n$ and $k_{\tau} = 1,...,9^{\tau -1} 8 .$ Moreover, $\{ T_n \}$ is decreasing and $\mathcal{L}^2 (T_n) = 9^{n} \ell^2 _{n+1} .$ We claim that
\begin{equation}\label{25.3.7_1}
\Gamma = \cap^{\infty} _{n=1} T_n .
\end{equation}
If this is true, we obtain from \eqref{25.2.5_1}  via $\sum^{\infty} _{n=1} n^{-2} = \pi^2 /6 $ that
\begin{equation*}
\mathcal{L}^2 (\Gamma) = \lim_{n \ra \infty} \mathcal{L}^2 (T_n) = \left[1- \left( \frac{\pi}{\ln 27}\right)^2 \right]^2>0 .
\end{equation*}
Obviously, every $\gamma^{k_{n}} _{n}$ is in $T_n$ for all $n \ge 1 .$ Take $z \in \mathcal{C}.$ If $z \notin \cap^{\infty} _{n=1} T_n , $ there is $n_0 \ge 1$ with $z \in (T_{n_0})^c .$ Hence, there is an open disk $B_0$ centered at  $z$ such that $ B_0 \cap T_{n} = \emptyset$ for all $n \ge n_0 .$ Therefore, $z$ is not an accumulation point of $\cup^{\infty} _{n=1} \cup^{9^{n-1} 8 } _{k_{n} =1 }  \gamma^{k_{n}} _{n} .$ This is impossible because $z \in \mathcal{C}.$ Thus $\Gamma \subset \cap^{\infty} _{n=1} T_n .$ 
To show the converse inclusion, we take $z \in \cap^{\infty} _{n=1} T_n \setminus (\cup^{\infty} _{n=1} \cup^{9^{n-1} 8 } _{k_{n} =1 }  \gamma^{k_{n}} _{n} ) .$ 
For every $n \ge 1, $ there is  $ Q^{j_n} _n \subset T_n$ such that $z \in Q^{j_n} _n $ and $Q^{j_{n+1}} _{n+1} \subset  Q^{j_n} _n .$
Denote by $z_n$ the middle point of a neighbouring $\gamma^{k_n} _n$ of $Q^{j_n} _n  .$
It follows from the geometry of $\{ T_n \}$ that 
\begin{equation*}
|z- z_n | \le \diam \, Q^{k_n} _n + d_n 
\qquad
\mbox{ and }
\qquad
|z_{n+1} - z_n| \le \diam \, Q^{k_n} _n + d_n .
\end{equation*}
Then $\{z_n\}$ is a Cauchy sequence and $ \lim_{n \ra \infty} z_n =z .$ So $z \in \mcc.$ We have that $\cap^{\infty} _{n=1} T_n \subset \Gamma .$ We conclude that  \eqref{25.3.7_1} holds.

\bigskip
{\it The pointwise growth conditions \eqref{25.2.5_3} and \eqref{25.2.7_1}.}
\bigskip

We will represent $Q \setminus \Gamma $ as a countable  union of   rectangles and crosses.
Let $D_n = [0, \ell_{n+1}) \times (-d_n /2,  d_n /2 ) .$  Then $D_n \cap \partial D_n$ is an open vertical segment of length $d_n .$
Define $C_n$ as the collection of points $(x,y) \in \mbr^2$ with $x \in [0, \ell_n)$
and 
\begin{equation*}
y \in 
\begin{cases}
(-\ell_{n+1} - \frac{d_n}{2} , \ell_{n+1} +\frac{d_n}{2} ), & \mbox{ if } x \in [2 \ell_{n+1} + d_n  , 2 \ell_{n+1} + 2 d_n] , \\
(-\frac{d_n}{2} ,  \frac{d_n}{2} ), & \mbox{ if } x \in [0, \ell_n) \setminus   [2 \ell_{n+1} + d_n  , 2 \ell_{n+1} + 2 d_n] .
\end{cases}
\end{equation*}
Then $C_n$ is cross-shaped and $C_n \cap \partial C_n $ is an open vertical segment.
Abusing  notation, we refer to the translated and rotated variants of  $C_n$ also by $C_n$ and do analogously for $D_n .$
To obtain $Q \setminus T_n $ for $n \ge 1 ,$
we remove $9^{n-1} 2$ suitable variants of  $D_n$ and the same number of variants of $C_n$  from $Q \setminus T_{n-1}.$ 
Denote by $\mathfrak{B}_n$ the collection of these $D_n $ and $C_n .$
It follows from \eqref{25.3.7_1} that
\begin{equation}\label{25.2.8_1}
Q \setminus \Gamma  = \cup^{\infty} _{n = 1}   \cup_{B \in \mathfrak{B}_n} B .
\end{equation}
We next give a decomposition of each $B_n \in \mathfrak{B}_n .$

Given  $B_n \in \mathfrak{B}_n ,$  we denote by $L_{n} = B_n \cap \partial B_n $ the open boundary segment. To be precise, $L_n = D_n \cap \partial D_n $ or $L_n = C_n \cap \partial  C_n  .$
For every $m > n , $ there is $B_m \in \mathfrak{B}_m$ such that its $L_m$-side is in $\partial B_n  .$
There are two such $B_m .$ Each such pair of $B_m$ is symmetric with respect to the middle axis of $B_n .$ 
For every such pair, we let $V_m$ be the subset of $B_n$ between their $L_m$ sides.
Denote by $\mathcal{V}_m$ the collection of all $V_m $ in $B_n .$
Then we decompose $B_n$ as 
\begin{equation}\label{25.2.13_1}
B_n = E_n  \cup \cup^{\infty} _{m=n+1} \cup_{V \in \mathcal{V}_m} V .
\end{equation}
We will apply this decomposition to prove \eqref{25.2.5_3}.

Consider  $z_0 :=-(1/2 ,0) \in \Omega.$ 
It suffices to prove \eqref{25.2.5_3} for this $z_0$ and all $z \in Q \cap \Omega .$ 
Notice that $Q \setminus \Gamma = Q \cap (\Omega \cup \Omega^* ) .$ 
By \eqref{25.2.8_1},  there is a $B_n \in \mathfrak{B}_n$ such that
$z \in B_n  .$
By the construction of $\Gamma ,$ there is a unique $B_{n_1} \in \mathfrak{B}_{n_1}$  with $1 \le n_1 <n$ such that $\overline{B_{n_1}} \cap B_{n} = L_n .$
Analogously, for $ B_{n_2} $ there is an unique $B_{n_2} \in \mathfrak{B}_{n_2}$ with $ 1 \le n_2 <n_1$ such that $\overline{B_{n_2}} \cap B_{n_1} = L_{n_1} .$ Continue this process until we obtain the $B_{n_{k_z}}$ with $n_{k_z} =1 $ and $\overline{B_{n_{k_z}}} \cap B_{n_{k_z -1} } = L_{n_{k_z -1}} .$
Let $z_n$ be the middle point of $L_n .$ 
Then
\begin{equation}\label{25.2.8_2}
h_{\Omega} (z_0 ,z) \le  
h_{\Omega} (z_0 ,z_1)
+ \sum^{k_z -1} _{i=1}  h_{\Omega} (z_{n_{i+1}} ,z_{n_i} )
+ h_{\Omega} (z_{n_1}, z_n)  
+ h_{\Omega} (z_{n} , z) .  
\end{equation}
We will establish an upper bound for this sum.

For $z \in B_n ,$ we denote by $z'$ the vertical projection  of $z$ in the middle axis of $B_n .$ For $V_m \in \mathcal{V}_m$ in \eqref{25.2.13_1}, we let $z_m$ be the midpoint of the $L_m$-side of $V_m .$ We associate a $z' _m$ with $z_m $ similarly as $z'$ was associated to $z.$
Let $z'' _m \in V_m$ with $\overrightarrow{z_m z'' _m } = \frac{d_m} {2|z_m- z' _m|} \overrightarrow{z_m z' _m}  .$ 
Denote by $U _m \subset V_m$ a rectangle whose one side is $L_m$ and the other side has length $|\overrightarrow{z_m z'' _m }| .$ 
Of course, for every $V_m ,$ there is a pair of $z_m,\ z''_m ,$ and $U_m .$ See Figure \ref{25.2.14_1} for these sets in $C_n .$
Depending on the location of $z \in B_n ,$ we next estimate $h_{\Omega} (z_{n} , z) .$

\begin{figure}
  \centering
  \includegraphics[width=\linewidth]{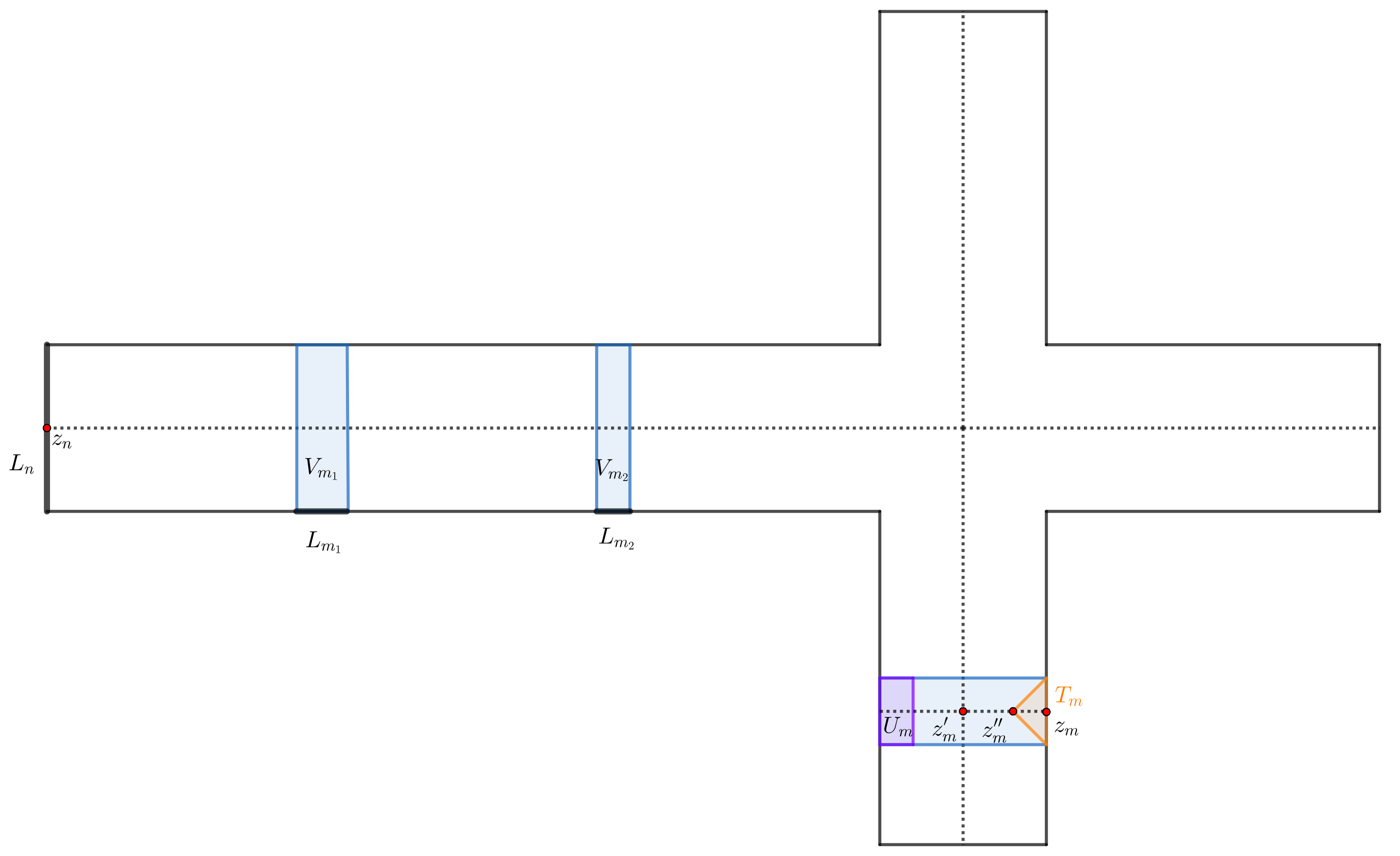}  
  \caption{A decomposition of $C_n$}
  \label{25.2.14_1}
\end{figure}

{\bf  Case 1: $z \in V_m \setminus U _m $.} We have
\begin{align*}
h_{\Omega} (z_{n} , z) \le &  h_{\Omega} (z_{n} , z'  ) + h_{\Omega} (z'  , z) \\
\le  &  \int_{\gamma_{z_n ,z' }} \frac{d \xi}{d(\xi ,\pao)} + \int_{L_{z'  ,z}} \frac{d \xi}{d(\xi ,\pao)} .
\end{align*}
Here $\gamma_{z_n ,z' } $ is the middle axis of $B_n $ between $z_n$ and $z' ,$ and $L_{z',z} \subset V_m$ is the segment joining $z'$ to $z .$
Notice that $d(\xi ,\pao) \approx  d(\xi , \partial B_n ) $ for $\xi \in \gamma_{z_n ,z' } \cup L_{z'  ,z} .$ 
Hence, by the values of $\ell_{n+1}$ and $d_n$ from  \eqref{25.2.5_1}
we have that 
\begin{align*}
\int_{\gamma_{z_n ,z' }} \frac{d \xi}{d(\xi ,\pao)} + \int_{L_{z'  ,z}} \frac{d \xi}{d(\xi ,\pao)} 
\lesssim  &  \frac{\ell_{n+1}}{d_n} + \log \frac{d_n }{2 d(z ,  \partial \Omega)}  \\
\lesssim & \log^2  \left( \frac{1 }{ d(z ,  \partial \Omega)} \right) . 
\end{align*}
By combining these two estimates, we arrive at 
$h_{\Omega} (z_{n} , z) \lesssim \log^2  \left( 1/ d(z ,  \partial \Omega) \right) .$

{\bf Case 2: $z \in U _m $.} 
By the result in Case 1, we obtain 
\begin{align}
h_{\Omega} (z_{n} , z) 
\le &  h_{\Omega} (z_{n} , z' _m ) + h_{\Omega} (z' _m , z'' _m) + h_{\Omega} ( z'' _m , z) \notag \\
\lesssim & \log^2  \left( \frac{1 }{ d(z'' _m ,  \partial \Omega)} \right) + h_{\Omega} ( z'' _m , z) \notag \\
\lesssim & \log^2  \left( \frac{1 }{ d(z ,  \partial \Omega)} \right) + h_{\Omega} ( z'' _m , z) , \label{25.3.10_1}
\end{align}
where the estimate $d(z ,\pao) \lesssim d(z'' _m ,\pao) $ was applied.
Denote by $T_m \subset U _m$ an isosceles right triangle with the hypotenuse $L_m$ and the right vertex $z'' _m .$
Given $z \in T_m ,$ we let $\tilde{z}$ be the projection of $z$ into $\overrightarrow{z_m z' _m } .$
Since $d(\xi ,\pao) \approx d(\xi , \partial L_m) $ for $\xi \in L_{z'' _m , \tilde{z}} \cup L_{\tilde{z} ,z} ,$
we obtain
\begin{align*}
h_{\Omega} ( z'' _m , z) 
\le & \int_{L_{z'' _m ,\tilde{z}}} \frac{d \xi}{d(\xi ,\pao)} +  \int_{L_{\tilde{z} ,z}} \frac{d \xi}{d(\xi ,\pao)} \\
\lesssim & \frac{|L_{z'' _m ,\tilde{z}}|}{d_m} + \log \frac{d_m}{2 d (z ,\pao)} \lesssim \log \frac{1}{ d (z ,\pao)}.
\end{align*} 
Analogously, we have $h_{\Omega} ( z'' _m , z) \lesssim \log 1/ d (z ,\pao)$ for $z \in U _m \setminus T_m .$ 
Together with \eqref{25.3.10_1}, we then obtain that $h_{\Omega} (z_{n} , z) \lesssim \log^2  (1/d(z ,  \partial \Omega) )$ for all $z \in U _m .$

{\bf Case 3: $z \in E_n$.} 
We have
\begin{equation*}
h_{\Omega} (z_n ,z) \le   h_{\Omega} (z_n ,z') + h_{\Omega} (z' ,z) 
\le  \int_{\gamma_{z_n ,z'}} \frac{d \xi}{d(\xi ,\pao)} + \int_{L_{z' ,z}}  \frac{d \xi}{d(\xi ,\pao)} .
\end{equation*}
Since $d(\xi ,\pao) \approx d(\xi ,\partial B_n)$ for $\xi \in \gamma_{z_n ,z'} \cup L_{z' ,z}  ,$ we have 
\begin{equation*}
\int_{\gamma_{z_n ,z'}} \frac{d \xi}{d(\xi ,\pao)} + \int_{L_{z' ,z}}  \frac{d \xi}{d(\xi ,\pao)} \lesssim \frac{\ell_{n+1}}{d_n} + \log \frac{d_n }{2 d(z ,\pao)} .
\end{equation*}
Therefore, $h_{\Omega} (z_n ,z) \lesssim \log^2  ( 1/ d(z ,  \partial \Omega) ) $ for all $z \in E_n .$

By  the estimates in the above three cases, we conclude that
\begin{equation}\label{25.2.13_5}
h_{\Omega} (z_{n} , z)  \le  \log^2  \left( \frac{1}{d(z ,  \partial \Omega)}  \right)
\end{equation} 
for all $z \in B_n .$

By arguments similar to those for Case 2, we also have  $h_{\Omega} (z_{n_{i+1}} ,z_{n_{i}}) \lesssim \log^2 ( 1/d_{n_{i+1}} ) .$
Therefore, 
\begin{align}
\sum^{k_z -1} _{i=1}  h_{\Omega} (z_{n_{i+1}} , z_{n_i} ) + h_{\Omega} (z_{n_1}, z_n) 
\lesssim &  \sum^n _{i=1} i^2 \lesssim n^3 \notag \\
\lesssim & \log^3 \left( \frac{1}{d(z ,\pao)} \right), \label{25.2.13_6}
\end{align}
for every $z \in B_n .$ 
It remains to estimate $h_{\Omega} (z_0 ,z_1)$ in \eqref{25.2.8_2}.
Since $\mathfrak{B}_1 \cap \Omega$ contains two elements, 
there are two choices for the location of $z_1 .$
In both of the case, one can show that $h_{\Omega} (z_0 ,z_1) \lesssim 1 $ by joining $z_0$ to $z_1$ by a suitable arc in $\Omega\setminus Q.$
By combining \eqref{25.2.8_2}, \eqref{25.2.13_5} with \eqref{25.2.13_6}, we then obtain
\begin{equation*}
h_{\Omega} (z_0 ,z) \le \log^3  \left( \frac{1}{d(z ,  \partial \Omega)}  \right)
\end{equation*}
for all $z \in B_n .$ This finishes the proof of \eqref{25.2.5_3}.

Fix $w_0 \in \Omega^* .$ 
Analogously to the proof of \eqref{25.2.5_3}, by the symmetry of our construction one can show that
\begin{equation*}
h_{\Omega^* } (w_0 , w) \lesssim \log ^3 \left( \frac{1}{d (w ,\pao)} \right)
\end{equation*}  
for all $w \in \Omega^* .$ This completes the proof of \eqref{25.2.7_1}.

Consider a given parametrization $\varphi$ of our Jordan curve embedded in $\mbs^2$ via $\mbs \subset \mbs^2 .$ By \eqref{25.2.5_3}, \eqref{25.2.7_1} and Theorem \ref{perus2} we may extend $\varphi$ to homeomorphisms $f^{+},\ f^{-}$ from the upper and lower parts of $\mbs^2 \setminus \mbs$ onto the complementary components of our Jordan curve, with $W^{1,p}$-regularity for all $p<2 .$ This gives us a self-homeomorphism $f : \mbs^2 \ra \mbs^2$ that belongs to $W^{1,p} (\mbs^2 \setminus \mbs , \mbs^2)$ for all $p<2 .$ It follows that $f \in W^{1,p} (\mbs^2  , \mbs^2) ,$ see e.g.  \cite{Jones 2000}.
\end{proof}

\bigskip
\bibliographystyle{amsplain}

\end{document}